\documentclass[11pt]{article}
\usepackage{amssymb,amsfonts,amsmath,graphicx,hyperref}

\newcommand{\mcm}[3]{\newcommand{#1}[#2]{{\ensuremath{#3}}}}

\usepackage{latexsym}
\usepackage{amssymb}

\mcm{\blank}{0}{(\emptybk)} \mcm{\dashbk}{0}{\mbox{---}}
\mcm{\emptybk}{0}{\:\:} \mcm{\hyph}{0}{\mbox{-}}

\mcm{\diagspace}{0}{\mbox{\hspace{2em}}}

\mcm{\cat}{1}{\mc{#1}} \mcm{\fcat}{1}{\mb{#1}}
\mcm{\mc}{1}{\mathcal{#1}} \mcm{\mr}{1}{\mathrm{#1}}
\mcm{\mi}{1}{\mathit{#1}} \mcm{\mb}{1}{\mathbf{#1}}
\mcm{\scat}{1}{\Bbb{#1}} \mcm{\twid}{1}{\widetilde{#1}}

\mcm{\elt}{0}{\in} \mcm{\sub}{0}{\,\subseteq\,}
\mcm{\such}{0}{\:|\:} \mcm{\without}{0}{\setminus}

\mcm{\atsr}{0}{\Box} \mcm{\eqv}{0}{\,\simeq\,}
\mcm{\iso}{0}{\,\cong\,}
\mcm{\of}{0}{\raisebox{0.2mm}{\ensuremath{\scriptstyle\circ}}}

\mcm{\bdry}{0}{\partial}

\mcm{\Bee}{0}{\cat{B}} \mcm{\Beep}{0}{\cat{B'}}
\mcm{\Eee}{0}{\cat{E}} \mcm{\Eeep}{0}{\cat{E'}}

\mcm{\Ess}{0}{\cat{S}} \mcm{\Tee}{0}{\cat{T}}
\mcm{\Teep}{0}{\cat{T'}} \mcm{\Stee}{0}{\scat{T}}
\mcm{\Steep}{0}{\scat{T'}}

\mcm{\blbk}{0}{\blank^{\blob}}
\mcm{\blob}{0}{\scriptscriptstyle{\bullet}}
\mcm{\stbk}{0}{\blank^{*}} \mcm{\ubl}{0}{{}^{\blob}}
\mcm{\ust}{0}{{}^{*}}

\mcm{\Cartpr}{0}{\pr{\Eee}{T}} \mcm{\Cartprp}{0}{\pr{\Eeep}{T'}}
\mcm{\Mnd}{0}{\triple{T}{\eta}{\mu}}
\mcm{\Zeropr}{0}{\pr{\Set}{\id}}

\mcm{\dopset}{0}{\ftrcat{\Delta^{\op}}{\Set}}
\mcm{\tropset}{0}{\ftrcat{\fcat{TR}^{\op}}{\Set}}

\mcm{\cod}{0}{\mr{cod}} \mcm{\dom}{0}{\mr{dom}}
\mcm{\End}{0}{\mr{End}} \mcm{\Hom}{0}{\mr{Hom}}
\mcm{\ob}{0}{\mr{ob}\,} \mcm{\op}{0}{\mr{op}}

\mcm{\comp}{0}{\mi{comp}} \mcm{\id}{0}{\mi{id}}
\mcm{\ids}{0}{\mi{ids}} \mcm{\mult}{0}{\mi{mult}}
\mcm{\unit}{0}{\mi{unit}}

\mcm{\Ab}{0}{\fcat{Ab}} \mcm{\Alg}{0}{\fcat{Alg}}
\mcm{\Bim}{1}{\fcat{Bim}(#1)} \mcm{\Cat}{0}{\fcat{Cat}}
\mcm{\Cay}{0}{\fcat{Cay}} \mcm{\Cpn}{1}{\pr{\Set/S_{#1}}{T_{#1}}}
\mcm{\fc}{0}{\fcat{fc}} \mcm{\fm}{0}{\fcat{fm}}
\mcm{\Graph}{0}{\fcat{Graph}} \mcm{\Gy}{0}{\fcat{Gy}}
\mcm{\Hpn}{1}{\pr{\Eee_{#1}}{P_{#1}}} \mcm{\Mon}{0}{\mb{Mon}}
\mcm{\Multicat}{0}{\fcat{Multicat}} \mcm{\One}{0}{\fcat{1}}
\mcm{\PD}{1}{\fcat{PD}_{#1}} \mcm{\Prof}{0}{\fcat{Prof}}
\mcm{\Set}{0}{\fcat{Set}} \mcm{\Span}{0}{\fcat{Span}}
\mcm{\Ssq}{0}{\fcat{Ssq}} \mcm{\Struc}{0}{\fcat{Struc}}
\mcm{\Sym}{0}{\fcat{Sym}} \mcm{\TR}{1}{\fcat{TR}(#1)}
\mcm{\Tr}{0}{\fcat{Tr}} \mcm{\Twocat}{0}{\fcat{2\hyph\Cat}}

\mcm{\integers}{0}{\mathbb{Z}}

\mcm{\range}{2}{#1,\,\ldots\,,#2}
\mcm{\bftuple}{2}{\tuplebts{\range{#1}{#2}}}
\mcm{\tuple}{3}{\tuplebts{\range{#1,#2}{#3}}}
\mcm{\rttuple}{1}{\tuplebts{\,\ldots\,,#1}}
\mcm{\abftuple}{2}{\atuplebts{\range{#1}{#2}}}
\mcm{\atuple}{3}{\atuplebts{\range{#1,#2}{#3}}}
\mcm{\arttuple}{1}{\atuplebts{\,\ldots\,,#1}}
\mcm{\sqbftuple}{2}{\obt\range{#1}{#2}\cbt}
\mcm{\pr}{2}{\tuplebts{#1,#2}}
\mcm{\triple}{3}{\tuplebts{#1,#2,#3}}

\mcm{\eend}{2}{#1[#2]} \mcm{\ehom}{3}{#1[#2,#3]}
\mcm{\ftrcat}{2}{[#1,#2]} \mcm{\homset}{3}{#1(#2,#3)}
\mcm{\multihom}{3}{#1(#2;#3)}
\mcm{\relhom}{5}{#1_{#2}(\range{#3}{#4};#5)}

\mcm{\go}{0}{\rTo} \mcm{\goby}{1}{\rTo^{#1}}
\mcm{\goesto}{0}{\,\longmapsto\,} \mcm{\goiso}{0}{\goby{\diso}}
\mcm{\monic}{0}{\rMonic} \mcm{\og}{0}{\lTo}
\mcm{\ogby}{1}{\lTo^{#1}}

\mcm{\gph}{2}{\spn{#1}{T #2}{#2}} \mcm{\graph}{4}{\spaan{#1}{T
#2}{#2}{#3}{#4}} \mcm{\oppair}{2}{\stackrel{\rTo^{#1}}{\lTo_{#2}}}
\mcm{\parpair}{2}{\stackrel{\rTo^{#1}}{\rTo_{#2}}}
\mcm{\spn}{3}{#2 \og #1 \go #3} \mcm{\spaan}{5}{#2 \ogby{#4} #1
\goby{#5} #3}

\mcm{\bktdvslob}{3}
    {\left(
    \begin{diagram}[height=1.5em]
    #1      \\
    \dTo>{\,#2} \\
    #3      \\
    \end{diagram}
    \right)}
\mcm{\slob}{3}{(#1 \goby{#2} #3)} \mcm{\vslob}{3}
    {\left.
    \begin{diagram}[height=1.5em]
    #1      \\
    \dTo>{\,#2} \\
    #3      \\
    \end{diagram}
    \right.}

\newenvironment{tree}
    {\begin{diagram}[height=1em,width=.75em,abut,noPS,tight]}
    {\end{diagram}}

\mcm{\enode}{0}{\circ}

\mcm{\nl}{1}{\stackrel{\textstyle #1}{\node}}
\mcm{\node}{0}{\bullet}

\mcm{\utree}{0}{\node}

\mcm{\diso}{0}{\sim}

\mcm{\vdiso}{0}{\wr}

\mcm{\nat}{0}{\mathbb{N}}

\mcm{\Onepr}{0}{\pr{\Graph}{\fc}}
\newlength{\nllwidth}
\newlength{\nllheight}
\newcommand{\stackbelow}[2]{%
\settowidth{\nllwidth}{\ensuremath{#1}\ensuremath{#2}}%
\settoheight{\nllheight}{\ensuremath{#2}}%
\addtolength{\nllheight}{.3ex}%
\mbox{%
\ensuremath{#1}%
\hspace{-.5\nllwidth}%
\raisebox{-1\nllheight}{\ensuremath{#2}}}}
\mcm{\nlal}{2}{\stackbelow{\nl{#1}}{#2}}
\mcm{\nll}{1}{\stackbelow{\node}{#1}} \mcm{\wun}{0}{\fcat{1}}
\mcm{\atuplebts}{1}{\langle #1 \rangle} \mcm{\tuplebts}{1}{(#1)}
\mcm{\bo}{0}{(} \mcm{\bc}{0}{)}
\mcm{\UBilax}{0}{\fcat{UBicat}_\mr{lax}}
\mcm{\UBiwk}{0}{\fcat{UBicat}_\mr{wk}}
\mcm{\UBistr}{0}{\fcat{UBicat}_\mr{str}}
\mcm{\Bilax}{0}{\fcat{Bicat}_\mr{lax}}
\mcm{\Biwk}{0}{\fcat{Bicat}_\mr{wk}}
\mcm{\Bistr}{0}{\fcat{Bicat}_\mr{str}} \mcm{\rotsub}{0}{\cup
\raisebox{0.1em}{$\scriptstyle{|}$}} \mcm{\pd}{0}{\fcat{pd}}
\mcm{\rep}{1}{\widehat{#1}} \mcm{\ovln}{1}{\overline{#1}}
\mcm{\Gph}{0}{\fcat{Gph}} \mcm{\tr}{0}{\fcat{tr}}

\mcm{\ladj}{0}{\,\dashv\,} \mcm{\zeropd}{0}{\node}
    {\end{diagram}}
\mcm{\END}{0}{\fcat{End}} \mcm{\HOM}{0}{\fcat{Hom}}

\newlength{\gwidth} 
\newlength{\gvert}  
\newlength{\gdrop}  
\newlength{\gbaredrop}  
\newlength{\goffset}    
\newlength{\gtemp}  

\newcommand{\present}[1]{%
\makebox[1\gwidth]{%
\rule[-1\gdrop]{0ex}{1\gvert}%
\raisebox{-1\gbaredrop}{#1}}}

\newcommand{\presentl}[1]{%
\makebox[1\gwidth][l]{%
\rule[-1\gdrop]{0ex}{1\gvert}%
\raisebox{-1\gbaredrop}{#1}}}

\newcommand{\presentr}[1]{%
\makebox[1\gwidth][r]{%
\rule[-1\gdrop]{0ex}{1\gvert}%
\raisebox{-1\gbaredrop}{#1}}}

\newcommand{\ginitdims}[2]{
\setlength{\unitlength}{1em}
\setlength{\goffset}{.25\unitlength}
\setlength{\gwidth}{#1\unitlength}
\setlength{\gvert}{#2\unitlength}
\setlength{\gdrop}{.5\gvert}
\addtolength{\gdrop}{-1\goffset}
\setlength{\gbaredrop}{1\gdrop}
\addtolength{\gvert}{.6\unitlength}
\addtolength{\gdrop}{.3\unitlength}}    

\newcommand{\cinitdims}[2]{
\setlength{\unitlength}{1em}
\setlength{\goffset}{.35\unitlength}
\setlength{\gwidth}{#1\unitlength}
\setlength{\gvert}{#2\unitlength}
\setlength{\gdrop}{.5\gvert}
\addtolength{\gdrop}{-1\goffset}
\setlength{\gbaredrop}{1\gdrop}
\addtolength{\gvert}{.6\unitlength}
\addtolength{\gdrop}{.3\unitlength}}    

\newcommand{\gsinitdims}[2]{
\setlength{\unitlength}{0.5em}
\setlength{\goffset}{.25\unitlength}
\setlength{\gwidth}{#1\unitlength}
\setlength{\gvert}{#2\unitlength}
\setlength{\gdrop}{.5\gvert}
\addtolength{\gdrop}{-1\goffset}
\setlength{\gbaredrop}{1\gdrop}
\addtolength{\gvert}{.6\unitlength}
\addtolength{\gdrop}{.3\unitlength}}    

\newcommand{\sidespic}[1]{%
\settowidth{\gtemp}{\ensuremath{#1}}%
\addtolength{\gwidth}{1\gtemp}}

\newcommand{\abovepic}[1]{%
\settoheight{\gtemp}{\ensuremath{#1}}%
\addtolength{\gvert}{1\gtemp}%
\settodepth{\gtemp}{\ensuremath{#1}}%
\addtolength{\gvert}{1\gtemp}}

\newcommand{\belowpic}[1]{%
\settoheight{\gtemp}{\ensuremath{#1}}%
\addtolength{\gvert}{1\gtemp}%
\addtolength{\gdrop}{1\gtemp}%
\settodepth{\gtemp}{\ensuremath{#1}}%
\addtolength{\gvert}{1\gtemp}%
\addtolength{\gdrop}{1\gtemp}}

\newcommand{\cell}[4]{\put(#1,#2){\makebox(0,0)[#3]{\ensuremath{#4}}}}
\mcm{\zmark}{0}{\scriptstyle{\bullet}}

\newcommand{\pregfst}[1]{%
\begin{picture}(0.5,0.2)(-0.5,-0.2)%
\cell{-0.1}{-0.2}{tr}{#1}%
\cell{0}{0}{c}{\zmark}%
\end{picture}}

\mcm{\gfst}{1}{%
\ginitdims{0.5}{0.4}%
\sidespic{#1}%
\belowpic{#1}%
\presentr{\pregfst{#1}}}

\newcommand{\preglst}[1]{%
\begin{picture}(0.5,0.2)(0,-0.2)%
\cell{0.1}{-0.2}{tl}{#1}%
\cell{0.05}{0}{c}{\zmark}%
\end{picture}}

\mcm{\glst}{1}{%
\ginitdims{.5}{.4}%
\sidespic{#1}%
\belowpic{#1}%
\presentl{\preglst{#1}}}

\newcommand{\preglft}[1]{%
\begin{picture}(0,0.2)(0,-0.2)%
\cell{-0.1}{-0.2}{tr}{#1}%
\cell{0.05}{0}{c}{\zmark}%
\end{picture}}

\mcm{\glft}{1}{%
\ginitdims{0}{.4}%
\belowpic{#1}%
\present{\preglft{#1}}}

\newcommand{\pregrgt}[1]{%
\begin{picture}(0,0.2)(0,-0.2)%
\cell{0.1}{-0.2}{tl}{#1}%
\cell{0.05}{0}{c}{\zmark}%
\end{picture}}

\mcm{\grgt}{1}{%
\ginitdims{0}{.4}%
\belowpic{#1}%
\present{\pregrgt{#1}}}

\newcommand{\pregblw}[1]{%
\begin{picture}(0,0.3)(0,-0.3)
\cell{0}{-0.3}{t}{#1}%
\cell{0.05}{0}{c}{\zmark}%
\end{picture}}

\mcm{\gblw}{1}{%
\ginitdims{0}{.6}%
\belowpic{#1}%
\present{\pregblw{#1}}}

\newcommand{\pregfbw}[1]{%
\begin{picture}(0,0.65)(0,-0.65)
\cell{0}{-0.65}{t}{#1}%
\cell{0.05}{0}{c}{\zmark}%
\end{picture}}

\mcm{\gfbw}{1}{%
\ginitdims{0}{1.3}%
\belowpic{#1}%
\present{\pregfbw{#1}}}

\newcommand{\pregzero}[1]{%
\begin{picture}(0.8,0.4)(-0.4,-0.4)
\cell{0}{-0.4}{t}{#1}%
\cell{0}{0}{c}{\zmark}%
\end{picture}}

\mcm{\gzero}{1}{%
\ginitdims{0.8}{.6}%
\belowpic{#1}%
\sidespic{#1}%
\present{\pregzero{#1}}}

\newcommand{\pregone}[1]{%
\begin{picture}(5,0.4)(0,-0.2)%
\cell{2.5}{0.2}{b}{#1}%
\put(0,0){\vector(1,0){5}}%
\end{picture}}

\mcm{\gone}{1}{%
\ginitdims{5}{0.4}%
\abovepic{#1}%
\present{\pregone{#1}}}

\newcommand{\pregtwo}[3]{%
\begin{picture}(5,3.4)(0,-0.2)%
\cell{2.5}{3.2}{b}{#1}%
\cell{2.5}{-.2}{t}{#2}%
\cell{2.7}{1.5}{l}{#3}%
\qbezier(0,1.5)(2.5,4.5)(5,1.5)%
\qbezier(0,1.5)(2.5,-1.5)(5,1.5)%
\put(5,1.5){\vector(1,-1){0}}%
\put(5,1.5){\vector(1,1){0}}%
\put(2.5,2.5){\vector(0,-1){2}}%
\end{picture}}

\mcm{\gtwo}{3}{%
\ginitdims{5}{3.4}%
\abovepic{#1}%
\belowpic{#2}%
\present{\pregtwo{#1}{#2}{#3}}}

\newcommand{\pregthree}[5]{%
\begin{picture}(5,5.4)(0,-1.2)%
\cell{2.5}{4.2}{b}{#1}%
\cell{1.5}{1.7}{b}{#2}%
\cell{2.5}{-1.2}{t}{#3}%
\cell{2.7}{2.75}{l}{#4}%
\cell{2.7}{0.25}{l}{#5}%
\qbezier(0,1.5)(2.5,6.5)(5,1.5)%
\qbezier(0,1.5)(2.5,-3.5)(5,1.5)%
\put(0,1.5){\vector(1,0){5}}%
\put(2.5,3.5){\vector(0,-1){1.5}}%
\put(2.5,1){\vector(0,-1){1.5}}%
\put(5,1.5){\vector(1,-3){0}}%
\put(5,1.5){\vector(1,3){0}}%
\end{picture}}

\mcm{\gthree}{5}{%
\ginitdims{5}{5.4}%
\abovepic{#1}%
\belowpic{#3}%
\present{\pregthree{#1}{#2}{#3}{#4}{#5}}}

\newcommand{\pregfour}[7]{%
\begin{picture}(5,8.4)(0,-2.7)%
\cell{2.5}{5.7}{b}{#1}%
\cell{1.5}{2.8}{b}{#2}%
\cell{1.5}{0.2}{t}{#3}%
\cell{2.5}{-2.7}{t}{#4}%
\cell{2.7}{4.25}{l}{#5}%
\cell{2.7}{1.5}{l}{#6}%
\cell{2.7}{-1.25}{l}{#7}%
\qbezier(0,1.5)(2.5,9.5)(5,1.5)%
\qbezier(0,1.5)(2.5,4)(5,1.5)%
\qbezier(0,1.5)(2.5,-1)(5,1.5)%
\qbezier(0,1.5)(2.5,-6.5)(5,1.5)%
\put(2.5,5.25){\vector(0,-1){2}}%
\put(2.5,2.5){\vector(0,-1){2}}%
\put(2.5,-0.25){\vector(0,-1){2}}%
\put(5,1.5){\vector(1,-4){0}}%
\put(5,1.5){\vector(4,-3){0}}%
\put(5,1.5){\vector(4,3){0}}%
\put(5,1.5){\vector(1,4){0}}%
\end{picture}}

\mcm{\gfour}{7}{%
\ginitdims{5}{8.4}%
\abovepic{#1}%
\belowpic{#4}%
\present{\pregfour{#1}{#2}{#3}{#4}{#5}{#6}{#7}}}

\newcommand{\pregthreecell}[5]{%
\begin{picture}(8,5)(-4,-2.5)%
\cell{0}{2.5}{b}{#1}%
\cell{0}{-2.5}{t}{#2}%
\cell{-1.7}{0}{r}{#3}%
\cell{1.7}{0}{l}{#4}%
\cell{0}{0.2}{b}{#5}%
\qbezier(-4,0)(0,4.2)(4,0)%
\qbezier(-4,0)(0,-4.2)(4,0)%
\qbezier(-0.5,1.8)(-2.5,0)(-0.5,-1.8)%
\qbezier(0.5,1.8)(2.5,0)(0.5,-1.8)%
\put(-1,0){\vector(1,0){2}}%
\put(4,0){\vector(1,-1){0}}%
\put(4,0){\vector(1,1){0}}%
\put(-0.5,-1.8){\vector(1,-1){0}}%
\put(0.5,-1.8){\vector(-1,-1){0}}%
\end{picture}}

\mcm{\gthreecell}{5}{%
\ginitdims{8}{5}%
\abovepic{#1}%
\belowpic{#2}%
\present{\pregthreecell{#1}{#2}{#3}{#4}{#5}}}

\newcommand{\pregthreecellu}{%
\begin{picture}(5,3.4)(-0.5,-0.2)%
\qbezier(-.5,1.5)(2,4.5)(4.5,1.5)%
\qbezier(-.5,1.5)(2,-1.5)(4.5,1.5)%
\qbezier(1.5,2.7)(0.5,1.5)(1.5,0.3)%
\qbezier(2.5,2.7)(3.5,1.5)(2.5,0.3)%
\put(1.3,1.5){\vector(1,0){1.4}}%
\put(4.5,1.5){\vector(1,-1){0}}%
\put(4.5,1.5){\vector(1,1){0}}%
\put(1.5,0.3){\vector(2,-3){0}}%
\put(2.5,0.3){\vector(-2,-3){0}}%
\end{picture}}

\mcm{\gthreecellu}{0}{%
\ginitdims{5}{3.4}%
\present{\pregthreecellu}}

\newcommand{\pregtwocentre}[3]{%
\begin{picture}(5,3.4)(0,-0.2)%
\cell{2.5}{3.2}{b}{#1}%
\cell{2.5}{-.2}{t}{#2}%
\cell{2.5}{1.5}{c}{#3}%
\qbezier(0,1.5)(2.5,4.5)(5,1.5)%
\qbezier(0,1.5)(2.5,-1.5)(5,1.5)%
\put(5,1.5){\vector(1,-1){0}}%
\put(5,1.5){\vector(1,1){0}}%
\put(2.5,2.5){\vector(0,-1){2}}%
\end{picture}}

\mcm{\gtwocentre}{3}{%
\ginitdims{5}{3.4}%
\abovepic{#1}%
\belowpic{#2}%
\present{\pregtwocentre{#1}{#2}{#3}}}

\newcommand{\pregspecialone}[9]{%
\begin{picture}(8,8)(-4,-4)%
\cell{0}{3.9}{b}{#1}%
\cell{-2}{-0.2}{t}{#2}%
\cell{0}{-3.9}{t}{#3}%
\cell{-1.5}{1.1}{r}{#4}%
\cell{0.2}{1.5}{l}{#5}%
\cell{1.5}{1.1}{l}{#6}%
\cell{0.2}{-2}{l}{#7}%
\cell{-0.9}{2.3}{b}{#8}%
\cell{0.9}{2.3}{b}{#9}%
\qbezier(-4,0)(0,8)(4,0)%
\qbezier(-4,0)(0,-8)(4,0)%
\qbezier(-0.5,3.4)(-3.5,2)(-0.5,0.6)%
\qbezier(0.5,3.4)(3.5,2)(0.5,0.6)%
\put(-4,0){\vector(1,0){8}}%
\put(0,3.4){\vector(0,-1){2.8}}%
\put(0,-0.8){\vector(0,-1){2.4}}%
\put(-1.5,2.2){\vector(1,0){1.2}}%
\put(0.3,2.2){\vector(1,0){1.2}}%
\put(4,0){\vector(1,-2){0}}%
\put(4,0){\vector(1,2){0}}%
\put(-0.5,0.6){\vector(2,-1){0}}%
\put(0.5,0.6){\vector(-2,-1){0}}%
\end{picture}}

\mcm{\gspecialone}{9}{%
\ginitdims{8}{8}%
\abovepic{#1}%
\belowpic{#3}%
\present{\pregspecialone{#1}{#2}{#3}{#4}{#5}{#6}{#7}{#8}{#9}}}

\newcommand{\pregspecialtwo}{%
\begin{picture}(5,3.4)(0,-0.2)%
\qbezier(0,1.5)(2.5,4.5)(5,1.5)%
\qbezier(0,1.5)(2.5,-1.5)(5,1.5)%
\qbezier(1.7,2.5)(0,1.5)(1.7,0.5)%
\qbezier(3.3,2.5)(5,1.5)(3.3,0.5)%
\put(5,1.5){\vector(1,-1){0}}%
\put(5,1.5){\vector(1,1){0}}%
\put(1.7,0.5){\vector(3,-2){0}}%
\put(3.3,0.5){\vector(-3,-2){0}}%
\put(2.5,2.5){\vector(0,-1){2}}%
\put(1.2,1.5){\vector(1,0){1}}%
\put(2.8,1.5){\vector(1,0){1}}%
\end{picture}}

\mcm{\gspecialtwo}{0}{%
\ginitdims{5}{3.4}%
\present{\pregspecialtwo}}

\newcommand{\pregspecialthree}{%
\begin{picture}(5,5.4)(0,-1.2)%
\qbezier(0,1.5)(2.5,6.5)(5,1.5)%
\qbezier(0,1.5)(2.5,-3.5)(5,1.5)%
\qbezier(2,3.5)(1,2.75)(2,2)%
\qbezier(3,3.5)(4,2.75)(3,2)%
\qbezier(2,1)(1,0.25)(2,-0.5)%
\qbezier(3,1)(4,0.25)(3,-0.5)%
\put(0,1.5){\vector(1,0){5}}%
\put(1.5,2.75){\vector(1,0){2}}%
\put(1.5,0.25){\vector(1,0){2}}%
\put(5,1.5){\vector(1,-3){0}}%
\put(5,1.5){\vector(1,3){0}}%
\put(2,2){\vector(1,-1){0}}%
\put(3,2){\vector(-1,-1){0}}%
\put(2,-0.5){\vector(1,-1){0}}%
\put(3,-0.5){\vector(-1,-1){0}}%
\end{picture}}

\mcm{\gspecialthree}{0}{%
\ginitdims{5}{5.4}%
\present{\pregspecialthree}}

\newcommand{\pregonew}[1]{%
\begin{picture}(8,0.4)(0,-0.2)%
\cell{4}{0.2}{b}{#1}%
\put(0,0){\vector(1,0){8}}%
\end{picture}}

\mcm{\gonew}{1}{%
\ginitdims{8}{0.4}%
\abovepic{#1}%
\present{\pregonew{#1}}}

\mcm{\gzersu}{0}{%
\gsinitdims{0}{.6}%
\present{\pregblw{}}}

\mcm{\gonesu}{0}{%
\gsinitdims{5}{0.4}%
\present{\pregone{}}}

\mcm{\gtwosu}{0}{%
\gsinitdims{5}{3.4}%
\present{\pregtwo{}{}{}}}

\mcm{\gthreesu}{0}{%
\gsinitdims{5}{5.4}%
\present{\pregthree{}{}{}{}{}}}

\mcm{\gfoursu}{0}{%
\gsinitdims{5}{8.4}%
\present{\pregfour{}{}{}{}{}{}{}}}

\newcommand{\precone}[1]{%
\begin{picture}(4.2,0.4)(-0.3,-0.2)%
\cell{1.8}{0.2}{b}{#1}%
\put(0,0){\vector(1,0){3.6}}%
\end{picture}}

\mcm{\cone}{1}{%
\cinitdims{4.2}{0.4}%
\abovepic{#1}%
\present{\precone{#1}}}

\mcm{\gfstsu}{0}{%
\gsinitdims{0.5}{0.4}%
\presentr{\pregfst{}}}

\mcm{\glstsu}{0}{%
\gsinitdims{0.5}{0.4}%
\presentl{\preglst{}}}

\newcommand{\prectwodbl}[3]%
{\begin{picture}(4.2,3.4)(-0.1,-0.2)%
\cell{2}{3.2}{b}{#1}%
\cell{2}{-0.2}{t}{#2}%
\cell{2.3}{1.5}{l}{#3}%
\qbezier(0,2)(2,4)(4,2)%
\qbezier(0,1)(2,-1)(4,1)%
\put(4,2){\vector(1,-1){0}}%
\put(4,1){\vector(1,1){0}}%
\put(1.9,2.5){\line(0,-1){1.8}}%
\put(2.1,2.5){\line(0,-1){1.8}}%
\cell{2.01}{0.4}{b}{\vee}%
\end{picture}}

\mcm{\ctwodbl}{3}{%
\cinitdims{4.2}{3.4}%
\abovepic{#1}%
\belowpic{#2}%
\present{\prectwodbl{#1}{#2}{#3}}}

\newcommand{\precthreedbl}[5]{%
\begin{picture}(4.2,5.4)(-0.1,-0.2)%
\cell{2}{5.2}{b}{#1}%
\cell{1}{2.7}{b}{#2}%
\cell{2}{-.2}{t}{#3}%
\cell{2.3}{3.75}{l}{#4}%
\cell{2.3}{1.25}{l}{#5}%
\qbezier(0,3)(2,7)(4,3)%
\qbezier(0,2)(2,-2)(4,2)%
\put(0,2.5){\vector(1,0){4}}%
\put(1.9,4.5){\line(0,-1){1.3}}%
\put(2.1,4.5){\line(0,-1){1.3}}%
\cell{2.01}{2.9}{b}{\vee}%
\put(1.9,2){\line(0,-1){1.3}}%
\put(2.1,2){\line(0,-1){1.3}}%
\cell{2.01}{0.4}{b}{\vee}%
\put(4,3){\vector(1,-3){0}}%
\put(4,2){\vector(1,3){0}}%
\end{picture}}

\mcm{\cthreedbl}{5}{%
\cinitdims{4.2}{5.4}%
\abovepic{#1}%
\belowpic{#3}%
\present{\precthreedbl{#1}{#2}{#3}{#4}{#5}}}

\newcommand{\precthreecelltrp}[5]{%
\begin{picture}(8.2,5)(-4.1,-2.5)%
\cell{0}{2.5}{b}{#1}%
\cell{0}{-2.5}{t}{#2}%
\cell{-1.8}{0}{r}{#3}%
\cell{1.8}{0}{l}{#4}%
\cell{0}{0.3}{b}{#5}%
\qbezier(-4,0.5)(0,4)(4,0.5)%
\qbezier(-4,-0.5)(0,-4)(4,-0.5)%
\qbezier(-0.6,2)(-2.6,0)(-0.6,-2)%
\qbezier(-0.4,2)(-2.4,0)(-0.5,-1.9)%
\cell{-0.6}{-2}{b}{\lrcorner}%
\qbezier(0.4,2)(2.4,0)(0.5,-1.9)%
\qbezier(0.6,2)(2.6,0)(0.6,-2)%
\cell{0.65}{-2}{b}{\llcorner}%
\put(-1,0.15){\line(1,0){1.7}}%
\put(-1,0){\line(1,0){2}}%
\put(-1,-0.15){\line(1,0){1.7}}%
\cell{1.15}{0}{r}{>}%
\put(4,0.5){\vector(1,-1){0}}%
\put(4,-0.5){\vector(1,1){0}}%
\end{picture}}

\mcm{\cthreecelltrp}{5}{%
\cinitdims{8.2}{5}%
\abovepic{#1}%
\belowpic{#2}%
\present{\precthreecelltrp{#1}{#2}{#3}{#4}{#5}}}

\newcommand{\prectwo}[3]%
{\begin{picture}(4.2,3.4)(-0.1,-0.2)%
\cell{2}{3.2}{b}{#1}%
\cell{2}{-0.2}{t}{#2}%
\cell{2.2}{1.5}{l}{#3}%
\qbezier(0,2)(2,4)(4,2)%
\qbezier(0,1)(2,-1)(4,1)%
\put(4,2){\vector(1,-1){0}}%
\put(4,1){\vector(1,1){0}}%
\put(2,2.5){\vector(0,-1){2}}%
\end{picture}}

\mcm{\ctwo}{3}{%
\cinitdims{4.2}{3.4}%
\abovepic{#1}%
\belowpic{#2}%
\present{\prectwo{#1}{#2}{#3}}}

\newcommand{\precthree}[5]{%
\begin{picture}(4.2,5.4)(-0.1,-0.2)%
\cell{2}{5.2}{b}{#1}%
\cell{1}{2.7}{b}{#2}%
\cell{2}{-.2}{t}{#3}%
\cell{2.2}{3.75}{l}{#4}%
\cell{2.2}{1.25}{l}{#5}%
\qbezier(0,3)(2,7)(4,3)%
\qbezier(0,2)(2,-2)(4,2)%
\put(0,2.5){\vector(1,0){4}}%
\put(2,4.5){\vector(0,-1){1.5}}%
\put(2,2){\vector(0,-1){1.5}}%
\put(4,3){\vector(1,-3){0}}%
\put(4,2){\vector(1,3){0}}%
\end{picture}}

\mcm{\cthree}{5}{%
\cinitdims{4.2}{5.4}%
\abovepic{#1}%
\belowpic{#3}%
\present{\precthree{#1}{#2}{#3}{#4}{#5}}}

\newcommand{\prectwoop}[3]%
{\begin{picture}(4.2,3.4)(-0.1,-0.2)%
\cell{2}{3.2}{b}{#1}%
\cell{2}{-0.2}{t}{#2}%
\cell{2.2}{1.5}{l}{#3}%
\qbezier(0,2)(2,4)(4,2)%
\qbezier(0,1)(2,-1)(4,1)%
\put(0,2){\vector(-1,-1){0}}%
\put(0,1){\vector(-1,1){0}}%
\put(2,2.5){\vector(0,-1){2}}%
\end{picture}}

\mcm{\ctwoop}{3}{%
\cinitdims{4.2}{3.4}%
\abovepic{#1}%
\belowpic{#2}%
\present{\prectwoop{#1}{#2}{#3}}}

\newcommand{\prectwopar}[4]{%
\begin{picture}(4.2,3.4)(-0.1,-0.2)%
\cell{2}{3.2}{b}{#1}%
\cell{2}{-0.2}{t}{#2}%
\cell{1.6}{1.5}{r}{#3}%
\cell{2.4}{1.5}{l}{#4}%
\qbezier(0,2)(2,4)(4,2)%
\qbezier(0,1)(2,-1)(4,1)%
\put(4,2){\vector(1,-1){0}}%
\put(4,1){\vector(1,1){0}}%
\put(1.8,2.5){\vector(0,-1){2}}%
\put(2.2,2.5){\vector(0,-1){2}}%
\end{picture}}

\mcm{\ctwopar}{4}{%
\cinitdims{4.2}{3.4}%
\abovepic{#1}%
\belowpic{#2}%
\present{\prectwopar{#1}{#2}{#3}{#4}}}

\newcommand{\precthreein}[5]{%
\begin{picture}(4.2,5.4)(-0.1,-0.2)%
\cell{2}{5.2}{b}{#1}%
\cell{1}{2.7}{b}{#2}%
\cell{2}{-.2}{t}{#3}%
\cell{2.2}{3.75}{l}{#4}%
\cell{2.2}{1.25}{l}{#5}%
\qbezier(0,3)(2,7)(4,3)%
\qbezier(0,2)(2,-2)(4,2)%
\put(0,2.5){\vector(1,0){4}}%
\put(2,4.5){\vector(0,-1){1.5}}%
\put(2,0.5){\vector(0,1){1.5}}%
\put(4,3){\vector(1,-3){0}}%
\put(4,2){\vector(1,3){0}}%
\end{picture}}

\mcm{\cthreein}{5}{%
\cinitdims{4.2}{5.4}%
\abovepic{#1}%
\belowpic{#3}%
\present{\precthreein{#1}{#2}{#3}{#4}{#5}}}

\newcommand{\precthreecell}[5]{%
\begin{picture}(8.2,5)(-4.1,-2.5)%
\cell{0}{2.5}{b}{#1}%
\cell{0}{-2.5}{t}{#2}%
\cell{-1.7}{0}{r}{#3}%
\cell{1.7}{0}{l}{#4}%
\cell{0}{0.2}{b}{#5}%
\qbezier(-4,0.5)(0,4)(4,0.5)%
\qbezier(-4,-0.5)(0,-4)(4,-0.5)%
\qbezier(-0.5,2)(-2.5,0)(-0.5,-2)%
\qbezier(0.5,2)(2.5,0)(0.5,-2)%
\put(-1,0){\vector(1,0){2}}%
\put(4,0.5){\vector(1,-1){0}}%
\put(4,-0.5){\vector(1,1){0}}%
\put(-0.5,-2){\vector(1,-1){0}}%
\put(0.5,-2){\vector(-1,-1){0}}%
\end{picture}}

\mcm{\cthreecell}{5}{%
\cinitdims{8.2}{5}%
\abovepic{#1}%
\belowpic{#2}%
\present{\precthreecell{#1}{#2}{#3}{#4}{#5}}}

\newcommand{\precthreecellpar}[6]{%
\begin{picture}(8.2,5)(-4.1,-2.5)%
\cell{0}{2.5}{b}{#1}%
\cell{0}{-2.5}{t}{#2}%
\cell{-1.7}{0}{r}{#3}%
\cell{1.7}{0}{l}{#4}%
\cell{0}{0.4}{b}{#5}%
\cell{0}{-0.4}{t}{#6}%
\qbezier(-4,0.5)(0,4)(4,0.5)%
\qbezier(-4,-0.5)(0,-4)(4,-0.5)%
\qbezier(-0.5,2)(-2.5,0)(-0.5,-2)%
\qbezier(0.5,2)(2.5,0)(0.5,-2)%
\put(-1,0.2){\vector(1,0){2}}%
\put(-1,-0.2){\vector(1,0){2}}%
\put(4,0.5){\vector(1,-1){0}}%
\put(4,-0.5){\vector(1,1){0}}%
\put(-0.5,-2){\vector(1,-1){0}}%
\put(0.5,-2){\vector(-1,-1){0}}%
\end{picture}}

\mcm{\cthreecellpar}{6}{%
\cinitdims{8.2}{5}%
\abovepic{#1}%
\belowpic{#2}%
\present{\precthreecellpar{#1}{#2}{#3}{#4}{#5}{#6}}}

\newcommand{\prectwov}[5]{%
\begin{picture}(3.4,4.2)(0.8,0.9)%
\cell{2.5}{5.1}{b}{#1}%
\cell{2.5}{0.9}{t}{#2}%
\cell{0.8}{3}{r}{#3}%
\cell{4.2}{3}{l}{#4}%
\cell{2.5}{3.2}{b}{#5}%
\qbezier(2,5)(0,3)(2,1)%
\qbezier(3,5)(5,3)(3,1)%
\put(2,1){\vector(1,-1){0}}%
\put(3,1){\vector(-1,-1){0}}%
\put(1.5,3){\vector(1,0){2}}%
\end{picture}}

\mcm{\ctwov}{5}{%
\cinitdims{3.4}{4.2}%
\abovepic{#1}%
\belowpic{#2}%
\sidespic{#3}%
\sidespic{#4}%
\present{\prectwov{#1}{#2}{#3}{#4}{#5}}}

\newcommand{\precthreecellv}[7]{%
\begin{picture}(5,8.2)(0.5,-1.6)%
\cell{3}{6.6}{b}{#1}%
\cell{3}{-1.6}{t}{#2}%
\cell{0.5}{2.5}{r}{#3}%
\cell{5.5}{2.5}{l}{#4}%
\cell{3}{4.2}{b}{#5}%
\cell{3}{0.8}{t}{#6}%
\cell{3.2}{2.5}{l}{#7}%
\qbezier(3.5,6.5)(7,2.5)(3.5,-1.5)%
\qbezier(2.5,6.5)(-1,2.5)(2.5,-1.5)%
\put(2.5,-1.5){\vector(1,-1){0}}%
\put(3.5,-1.5){\vector(-1,-1){0}}%
\qbezier(1,3)(3,5)(5,3)%
\qbezier(1,2)(3,0)(5,2)%
\put(5,3){\vector(1,-1){0}}%
\put(5,2){\vector(1,1){0}}%
\put(3,3.5){\vector(0,-1){2}}%
\end{picture}}

\mcm{\cthreecellv}{7}{%
\cinitdims{5}{8.2}%
\abovepic{#1}%
\belowpic{#2}%
\sidespic{#3}%
\sidespic{#4}%
\present{\precthreecellv{#1}{#2}{#3}{#4}{#5}{#6}{#7}}}

\newcommand{\pretopez}[2]{%
\begin{picture}(2.6,2.3)(-1.3,-2.2)%
\cell{0}{-2.2}{t}{#1}%
\cell{0}{-1.2}{c}{#2}%
\qbezier(0,0)(-2,-2)(0,-2)%
\qbezier(0,0)(2,-2)(0,-2)%
\put(0,0){\vector(-1,1){0}}%
\end{picture}}

\mcm{\topez}{2}{%
\ginitdims{2.6}{2.3}%
\belowpic{#1}%
\present{\pretopez{#1}{#2}}}

\newcommand{\pretopea}[3]{%
\begin{picture}(4,1.9)(-2,-0,2)%
\cell{0}{1.7}{b}{#1}%
\cell{0}{-0.2}{t}{#2}%
\cell{0}{0.7}{c}{#3}%
\qbezier(-2,0)(0,3)(2,0)%
\put(-2,0){\vector(1,0){4}}%
\put(2,0){\vector(2,-3){0}}%
\end{picture}}

\mcm{\topea}{3}{%
\ginitdims{4}{1.9}%
\abovepic{#1}%
\belowpic{#2}%
\present{\pretopea{#1}{#2}{#3}}}

\newcommand{\pretopeb}[4]{%
\begin{picture}(4,2.2)(-2,-0.2)%
\cell{-1.1}{1}{br}{#1}%
\cell{1.1}{1}{bl}{#2}%
\cell{0}{-0.2}{t}{#3}%
\cell{0}{0.8}{c}{#4}%
\put(-2,0){\vector(1,1){2}}%
\put(0,2){\vector(1,-1){2}}%
\put(-2,0){\vector(1,0){4}}%
\end{picture}}

\mcm{\topeb}{4}{%
\ginitdims{4}{2.2}%
\belowpic{#3}%
\present{\pretopeb{#1}{#2}{#3}{#4}}}

\newcommand{\pretopec}[5]{%
\begin{picture}(4,2.2)(-2,-0.2)%
\cell{-1.8}{1}{br}{#1}%
\cell{0}{2.2}{b}{#2}%
\cell{1.8}{1}{bl}{#3}%
\cell{0}{-0.2}{t}{#4}%
\cell{0}{0.8}{c}{#5}%
\put(-2,0){\vector(1,2){1}}%
\put(-1,2){\vector(1,0){2}}%
\put(1,2){\vector(1,-2){1}}%
\put(-2,0){\vector(1,0){4}}%
\end{picture}}

\mcm{\topec}{5}{%
\ginitdims{4}{2.2}%
\sidespic{#1}%
\abovepic{#2}%
\sidespic{#3}%
\belowpic{#4}%
\present{\pretopec{#1}{#2}{#3}{#4}{#5}}}

\newcommand{\pretoped}[6]{%
\begin{picture}(4,2.5)(-2,-0.2)%
\cell{-2}{0.6}{br}{#1}%
\cell{-0.7}{2.2}{br}{#2}%
\cell{0.7}{2.2}{bl}{#3}%
\cell{2}{0.6}{bl}{#4}%
\cell{0}{-0.2}{t}{#5}%
\cell{0}{0.8}{c}{#6}%
\put(-2,0){\vector(1,3){0.5}}%
\put(-1.5,1.5){\vector(3,2){1.5}}%
\put(0,2.5){\vector(3,-2){1.5}}%
\put(1.5,1.5){\vector(1,-3){0.5}}%
\put(-2,0){\vector(1,0){4}}%
\end{picture}}

\mcm{\toped}{6}{%
\ginitdims{4}{2.5}%
\sidespic{#1}%
\abovepic{#2}%
\abovepic{#3}%
\sidespic{#4}%
\belowpic{#5}%
\present{\pretoped{#1}{#2}{#3}{#4}{#5}{#6}}}

\newcommand{\pretopeq}[5]{%
\begin{picture}(4,2.5)(-2,-0.2)%
\cell{-2}{0.6}{br}{#1}%
\cell{-1}{2.2}{br}{#2}%
\cell{2}{0.6}{bl}{#3}%
\cell{0}{-0.2}{t}{#4}%
\cell{0}{0.8}{c}{#5}%
\put(-2,0){\vector(1,3){0.5}}%
\put(-1.5,1.5){\vector(1,1){1}}%
\cell{0.9}{2.3}{c}{\ddots}
\put(1.5,1.5){\vector(1,-3){0.5}}%
\put(-2,0){\vector(1,0){4}}%
\end{picture}}

\mcm{\topeq}{5}{%
\ginitdims{4}{2.5}%
\sidespic{#1}%
\abovepic{#2}%
\sidespic{#3}%
\belowpic{#4}%
\present{\pretopeq{#1}{#2}{#3}{#4}{#5}}}

\newcommand{\pretopebase}[1]{%
\begin{picture}(4,0.4)(0,-0.2)%
\cell{2}{0.2}{b}{#1}%
\put(0,0){\vector(1,0){4}}%
\end{picture}}

\mcm{\topebase}{1}{%
\ginitdims{4}{0.4}%
\abovepic{#1}%
\present{\pretopebase{#1}}}

\newcommand{\pretopezs}[2]{%
\begin{picture}(2.6,2.3)(-1.3,-2.2)%
\cell{0}{-2.2}{t}{#1}%
\cell{0}{-1.2}{c}{#2}%
\qbezier(0,0)(-2,-2)(0,-2)%
\qbezier(0,0)(2,-2)(0,-2)%
\end{picture}}

\mcm{\topezs}{2}{%
\ginitdims{2.6}{2.3}%
\belowpic{#1}%
\present{\pretopezs{#1}{#2}}}

\newcommand{\pretopeas}[3]{%
\begin{picture}(4,1.9)(-2,-0,2)%
\cell{0}{1.7}{b}{#1}%
\cell{0}{-0.2}{t}{#2}%
\cell{0}{0.7}{c}{#3}%
\qbezier(-2,0)(0,3)(2,0)%
\put(-2,0){\line(1,0){4}}%
\end{picture}}

\mcm{\topeas}{3}{%
\ginitdims{4}{1.9}%
\abovepic{#1}%
\belowpic{#2}%
\present{\pretopeas{#1}{#2}{#3}}}

\newcommand{\pretopebs}[4]{%
\begin{picture}(4,2.2)(-2,-0.2)%
\cell{-1.1}{1}{br}{#1}%
\cell{1.1}{1}{bl}{#2}%
\cell{0}{-0.2}{t}{#3}%
\cell{0}{0.8}{c}{#4}%
\put(-2,0){\line(1,1){2}}%
\put(0,2){\line(1,-1){2}}%
\put(-2,0){\line(1,0){4}}%
\end{picture}}

\mcm{\topebs}{4}{%
\ginitdims{4}{2.2}%
\belowpic{#3}%
\present{\pretopebs{#1}{#2}{#3}{#4}}}

\newcommand{\pretopecs}[5]{%
\begin{picture}(4,2.2)(-2,-0.2)%
\cell{-1.8}{1}{br}{#1}%
\cell{0}{2.2}{b}{#2}%
\cell{1.8}{1}{bl}{#3}%
\cell{0}{-0.2}{t}{#4}%
\cell{0}{0.8}{c}{#5}%
\put(-2,0){\line(1,2){1}}%
\put(-1,2){\line(1,0){2}}%
\put(1,2){\line(1,-2){1}}%
\put(-2,0){\line(1,0){4}}%
\end{picture}}

\mcm{\topecs}{5}{%
\ginitdims{4}{2.2}%
\sidespic{#1}%
\abovepic{#2}%
\sidespic{#3}%
\belowpic{#4}%
\present{\pretopecs{#1}{#2}{#3}{#4}{#5}}}

\newcommand{\pretopeds}[6]{%
\begin{picture}(4,2.5)(-2,-0.2)%
\cell{-2}{0.6}{br}{#1}%
\cell{-0.7}{2.2}{br}{#2}%
\cell{0.7}{2.2}{bl}{#3}%
\cell{2}{0.6}{bl}{#4}%
\cell{0}{-0.2}{t}{#5}%
\cell{0}{0.8}{c}{#6}%
\put(-2,0){\line(1,3){0.5}}%
\put(-1.5,1.5){\line(3,2){1.5}}%
\put(0,2.5){\line(3,-2){1.5}}%
\put(1.5,1.5){\line(1,-3){0.5}}%
\put(-2,0){\line(1,0){4}}%
\end{picture}}

\mcm{\topeds}{6}{%
\ginitdims{4}{2.5}%
\sidespic{#1}%
\abovepic{#2}%
\abovepic{#3}%
\sidespic{#4}%
\belowpic{#5}%
\present{\pretopeds{#1}{#2}{#3}{#4}{#5}{#6}}}

\newcommand{\pretopeqs}[5]{%
\begin{picture}(4,2.5)(-2,-0.2)%
\cell{-2}{0.6}{br}{#1}%
\cell{-1}{2.2}{br}{#2}%
\cell{2}{0.6}{bl}{#3}%
\cell{0}{-0.2}{t}{#4}%
\cell{0}{0.8}{c}{#5}%
\put(-2,0){\line(1,3){0.5}}%
\put(-1.5,1.5){\line(1,1){1}}%
\cell{0.9}{2.3}{c}{\ddots}
\put(1.5,1.5){\line(1,-3){0.5}}%
\put(-2,0){\line(1,0){4}}%
\end{picture}}

\mcm{\topeqs}{5}{%
\ginitdims{4}{2.5}%
\sidespic{#1}%
\abovepic{#2}%
\sidespic{#3}%
\belowpic{#4}%
\present{\pretopeqs{#1}{#2}{#3}{#4}{#5}}}

\newcommand{\pretopebases}[1]{%
\begin{picture}(4,0.4)(0,-0.2)%
\cell{2}{0.2}{b}{#1}%
\put(0,0){\line(1,0){4}}%
\end{picture}}

\mcm{\topebases}{1}{%
\ginitdims{4}{0.4}%
\abovepic{#1}%
\present{\pretopebases{#1}}}

\newcommand{\pregdots}[6]{%
\begin{picture}(5,8.4)(0,-2.7)%
\cell{2.5}{5.7}{b}{#1}%
\cell{1.5}{2.8}{b}{#2}%
\cell{1.5}{0.2}{t}{#3}%
\cell{2.5}{-2.7}{t}{#4}%
\cell{2.7}{4.25}{l}{#5}%
\cell{2.7}{-1.25}{l}{#6}%
\qbezier(0,1.5)(2.5,9.5)(5,1.5)%
\qbezier(0,1.5)(2.5,4)(5,1.5)%
\qbezier(0,1.5)(2.5,-1)(5,1.5)%
\qbezier(0,1.5)(2.5,-6.5)(5,1.5)%
\put(2.5,5.25){\vector(0,-1){2}}%
\put(2.5,-0.25){\vector(0,-1){2}}%
\cell{2.5}{1.7}{c}{\vdots}%
\put(5,1.5){\vector(1,-4){0}}%
\put(5,1.5){\vector(4,-3){0}}%
\put(5,1.5){\vector(4,3){0}}%
\put(5,1.5){\vector(1,4){0}}%
\end{picture}}

\mcm{\gdots}{6}{%
\ginitdims{5}{8.4}%
\abovepic{#1}%
\belowpic{#4}%
\present{\pregdots{#1}{#2}{#3}{#4}{#5}{#6}}}

\newlength{\volt}
\makeatletter

\def\diagram{\m@th\leftwidth=\z@ \rightwidth=\z@ \topheight=\z@
\botheight=\z@ \setbox\@picbox\hbox\bgroup}

\def\enddiagram{\egroup\wd\@picbox\rightwidth\unitlength
\ht\@picbox\topheight\unitlength \dp\@picbox\botheight\unitlength
\hskip\leftwidth\unitlength\box\@picbox}

\def\bfig{\begin{diagram}}
\def\efig{\end{diagram}}
\newcount\wideness \newcount\leftwidth \newcount\rightwidth
\newcount\highness \newcount\topheight \newcount\botheight

\def\ratchet#1#2{\ifnum#1<#2 \global #1=#2 \fi}

\def\putbox(#1,#2)#3{%
\horsize{\wideness}{#3} \divide\wideness by 2 {\advance\wideness
by #1 \ratchet{\rightwidth}{\wideness}} {\advance\wideness by -#1
\ratchet{\leftwidth}{\wideness}} \vertsize{\highness}{#3}
\divide\highness by 2 {\advance\highness by #2
\ratchet{\topheight}{\highness}} {\advance\highness by -#2
\ratchet{\botheight}{\highness}} \put(#1,#2){\makebox(0,0){$#3$}}}

\def\putlbox(#1,#2)#3{%
\horsize{\wideness}{#3} {\advance\wideness by #1
\ratchet{\rightwidth}{\wideness}} {\ratchet{\leftwidth}{-#1}}
\vertsize{\highness}{#3} \divide\highness by 2 {\advance\highness
by #2 \ratchet{\topheight}{\highness}} {\advance\highness by -#2
\ratchet{\botheight}{\highness}}
\put(#1,#2){\makebox(0,0)[l]{$#3$}}}

\def\putrbox(#1,#2)#3{%
\horsize{\wideness}{#3} {\ratchet{\rightwidth}{#1}}
{\advance\wideness by -#1 \ratchet{\leftwidth}{\wideness}}
\vertsize{\highness}{#3} \divide\highness by 2 {\advance\highness
by #2 \ratchet{\topheight}{\highness}} {\advance\highness by -#2
\ratchet{\botheight}{\highness}}
\put(#1,#2){\makebox(0,0)[r]{$#3$}}}

\def\adjust[#1]{} 

\newcount \coefa
\newcount \coefb
\newcount \coefc
\newcount\tempcounta
\newcount\tempcountb
\newcount\tempcountc
\newcount\tempcountd
\newcount\xext
\newcount\yext
\newcount\xoff
\newcount\yoff
\newcount\gap%
\newcount\arrowtypea
\newcount\arrowtypeb
\newcount\arrowtypec
\newcount\arrowtyped
\newcount\arrowtypee
\newcount\height
\newcount\width
\newcount\xpos
\newcount\ypos
\newcount\run
\newcount\rise
\newcount\arrowlength
\newcount\halflength
\newcount\arrowtype
\newdimen\tempdimen
\newdimen\xlen
\newdimen\ylen
\newsavebox{\tempboxa}%
\newsavebox{\tempboxb}%
\newsavebox{\tempboxc}%

\newdimen\w@dth

\def\setw@dth#1#2{\setbox\z@\hbox{\m@th$#1$}\w@dth=\wd\z@
\setbox\@ne\hbox{\m@th$#2$}\ifnum\w@dth<\wd\@ne \w@dth=\wd\@ne \fi
\advance\w@dth by 1.2em}

\def\t@^#1_#2{\allowbreak\def\n@one{#1}\def\n@two{#2}\mathrel
{\setw@dth{#1}{#2} \mathop{\hbox to
\w@dth{\rightarrowfill}}\limits \ifx\n@one\empty\else
^{\box\z@}\fi \ifx\n@two\empty\else _{\box\@ne}\fi}}
\def\t@@^#1{\@ifnextchar_{\t@^{#1}}{\t@^{#1}_{}}}
\def\to{\@ifnextchar^{\t@@}{\t@@^{}}}

\def\t@left^#1_#2{\def\n@one{#1}\def\n@two{#2}\mathrel{\setw@dth{#1}{#2}
\mathop{\hbox to \w@dth{\leftarrowfill}}\limits
\ifx\n@one\empty\else ^{\box\z@}\fi \ifx\n@two\empty\else
_{\box\@ne}\fi}}
\def\t@@left^#1{\@ifnextchar_{\t@left^{#1}}{\t@left^{#1}_{}}}
\def\toleft{\@ifnextchar^{\t@@left}{\t@@left^{}}}

\def\two@^#1_#2{\allowbreak
\def\n@one{#1}\def\n@two{#2}\mathrel{\setw@dth{#1}{#2}
\mathop{\vcenter{\lineskip\z@\baselineskip\z@
                 \hbox to \w@dth{\rightarrowfill}%
                 \hbox to \w@dth{\rightarrowfill}}%
       }\limits
\ifx\n@one\empty\else ^{\box\z@}\fi \ifx\n@two\empty\else
_{\box\@ne}\fi}}
\def\tw@@^#1{\@ifnextchar _{\two@^{#1}}{\two@^{#1}_{}}}
\def\two{\@ifnextchar ^{\tw@@}{\tw@@^{}}}

\def\tofr@^#1_#2{\def\n@one{#1}\def\n@two{#2}\mathrel{\setw@dth{#1}{#2}
\mathop{\vcenter{\hbox to \w@dth{\rightarrowfill}\kern-1.7ex
                 \hbox to \w@dth{\leftarrowfill}}%
       }\limits
\ifx\n@one\empty\else ^{\box\z@}\fi \ifx\n@two\empty\else
_{\box\@ne}\fi}}
\def\t@fr@^#1{\@ifnextchar_ {\tofr@^{#1}}{\tofr@^{#1}_{}}}
\def\tofro{\@ifnextchar^ {\t@fr@}{\t@fr@^{}}}

\def\mon{\mathop{\m@th\hbox to
      14.6\P@{\lasyb\char'51\hskip-2.1\P@$\arrext$\hss
$\mathord\rightarrow$}}\limits} 
\def\leftmono{\mathrel{\m@th\hbox to
14.6\P@{$\mathord\leftarrow$\hss$\arrext$\hskip-2.1\P@\lasyb\char'50%
}}\limits} 
\mathchardef\arrext="0200       

\def\settypes(#1,#2,#3){\arrowtypea#1 \arrowtypeb#2 \arrowtypec#3}
\def\settoheight#1#2{\setbox\@tempboxa\hbox{#2}#1\ht\@tempboxa\relax}%
\def\settodepth#1#2{\setbox\@tempboxa\hbox{#2}#1\dp\@tempboxa\relax}%
\def\settokens`#1`#2`#3`#4`{%
     \def\tokena{#1}\def\tokenb{#2}\def\tokenc{#3}\def\tokend{#4}}
\def\setsqparms[#1`#2`#3`#4;#5`#6]{%
\arrowtypea #1 \arrowtypeb #2 \arrowtypec #3 \arrowtyped #4
\width #5 \height #6 }
\def\setpos(#1,#2){\xpos=#1 \ypos#2}

\def\settriparms[#1`#2`#3;#4]{\settripairparms[#1`#2`#3`1`1;#4]}%

\def\settripairparms[#1`#2`#3`#4`#5;#6]{%
\arrowtypea #1 \arrowtypeb #2 \arrowtypec #3 \arrowtyped #4
\arrowtypee #5 \width #6 \height #6 }

\def\resetparms{\settripairparms[1`1`1`1`1;500]\width 500}

\resetparms

\def\mvector(#1,#2)#3{
\put(0,0){\vector(#1,#2){#3}}%
\put(0,0){\vector(#1,#2){26}}%
}
\def\evector(#1,#2)#3{{
\arrowlength #3
\put(0,0){\vector(#1,#2){\arrowlength}}%
\advance \arrowlength by-30
\put(0,0){\vector(#1,#2){\arrowlength}}%
}}

\def\horsize#1#2{%
\settowidth{\tempdimen}{$#2$}%
#1=\tempdimen \divide #1 by\unitlength }

\def\vertsize#1#2{%
\settoheight{\tempdimen}{$#2$}%
#1=\tempdimen
\settodepth{\tempdimen}{$#2$}%
\advance #1 by\tempdimen \divide #1 by\unitlength }

\def\putvector(#1,#2)(#3,#4)#5#6{{%
\ifnum3<\arrowtype \putdashvector(#1,#2)(#3,#4)#5\arrowtype \else
\ifnum\arrowtype<-3 \putdashvector(#1,#2)(#3,#4)#5\arrowtype \else
\xpos=#1 \ypos=#2 \run=#3 \rise=#4 \arrowlength=#5 \ifnum
\arrowtype<0
    \ifnum \run=0
        \advance \ypos by-\arrowlength
    \else
        \tempcounta \arrowlength
        \multiply \tempcounta by\rise
        \divide \tempcounta by\run
        \ifnum\run>0
            \advance \xpos by\arrowlength
            \advance \ypos by\tempcounta
        \else
            \advance \xpos by-\arrowlength
            \advance \ypos by-\tempcounta
        \fi
    \fi
    \multiply \arrowtype by-1
    \multiply \rise by-1
    \multiply \run by-1
\fi \ifcase \arrowtype
\or \put(\xpos,\ypos){\vector(\run,\rise){\arrowlength}}%
\or \put(\xpos,\ypos){\mvector(\run,\rise)\arrowlength}%
\or \put(\xpos,\ypos){\evector(\run,\rise){\arrowlength}}%
\fi\fi\fi }}

\def\putsplitvector(#1,#2)#3#4{
\xpos #1 \ypos #2 \arrowtype #4 \halflength #3 \arrowlength #3
\gap 140 \advance \halflength by-\gap \divide \halflength by2
\ifnum\arrowtype>0
   \ifcase \arrowtype
   \or \put(\xpos,\ypos){\line(0,-1){\halflength}}%
       \advance\ypos by-\halflength
       \advance\ypos by-\gap
       \put(\xpos,\ypos){\vector(0,-1){\halflength}}%
   \or \put(\xpos,\ypos){\line(0,-1)\halflength}%
       \put(\xpos,\ypos){\vector(0,-1)3}%
       \advance\ypos by-\halflength
       \advance\ypos by-\gap
       \put(\xpos,\ypos){\vector(0,-1){\halflength}}%
   \or \put(\xpos,\ypos){\line(0,-1)\halflength}%
       \advance\ypos by-\halflength
       \advance\ypos by-\gap
       \put(\xpos,\ypos){\evector(0,-1){\halflength}}%
   \fi
\else \arrowtype=-\arrowtype
   \ifcase\arrowtype
   \or \advance \ypos by-\arrowlength
       \put(\xpos,\ypos){\line(0,1){\halflength}}%
       \advance\ypos by\halflength
       \advance\ypos by\gap
       \put(\xpos,\ypos){\vector(0,1){\halflength}}%
   \or \advance \ypos by-\arrowlength
       \put(\xpos,\ypos){\line(0,1)\halflength}%
       \put(\xpos,\ypos){\vector(0,1)3}%
       \advance\ypos by\halflength
       \advance\ypos by\gap
       \put(\xpos,\ypos){\vector(0,1){\halflength}}%
   \or \advance \ypos by-\arrowlength
       \put(\xpos,\ypos){\line(0,1)\halflength}%
       \advance\ypos by\halflength
       \advance\ypos by\gap
       \put(\xpos,\ypos){\evector(0,1){\halflength}}%
   \fi
\fi }

\def\putmorphism(#1)(#2,#3)[#4`#5`#6]#7#8#9{{%
\run #2 \rise #3 \ifnum\rise=0
  \puthmorphism(#1)[#4`#5`#6]{#7}{#8}#9%
\else\ifnum\run=0
  \putvmorphism(#1)[#4`#5`#6]{#7}{#8}#9%
\else
\setpos(#1)%
\arrowlength #7 \arrowtype #8 \ifnum\run=0 \else\ifnum\rise=0
\else \ifnum\run>0
    \coefa=1
\else
   \coefa=-1
\fi \ifnum\arrowtype>0
   \coefb=0
   \coefc=-1
\else
   \coefb=\coefa
   \coefc=1
   \arrowtype=-\arrowtype
\fi \width=2 \multiply \width by\run \divide \width by\rise
\ifnum \width<0  \width=-\width\fi \advance\width by60 \if l#9
\width=-\width\fi
\putbox(\xpos,\ypos){#4}
{\multiply \coefa by\arrowlength
\advance\xpos by\coefa \multiply \coefa by\rise \divide \coefa
by\run \advance \ypos by\coefa
\putbox(\xpos,\ypos){#5} }%
{\multiply \coefa by\arrowlength
\divide \coefa by2 \advance \xpos by\coefa \advance \xpos by\width
\multiply \coefa by\rise \divide \coefa by\run \advance \ypos
by\coefa
\if l#9%
   \putrbox(\xpos,\ypos){#6}%
\else\if r#9%
   \putlbox(\xpos,\ypos){#6}%
\fi\fi }%
{\multiply \rise by-\coefc
\multiply \run by-\coefc \multiply \coefb by\arrowlength \advance
\xpos by\coefb \multiply \coefb by\rise \divide \coefb by\run
\advance \ypos by\coefb \multiply \coefc by70 \advance \ypos
by\coefc \multiply \coefc by\run \divide \coefc by\rise \advance
\xpos by\coefc \multiply \coefa by140 \multiply \coefa by\run
\divide \coefa by\rise \advance \arrowlength by\coefa
\ifcase\arrowtype
\or \put(\xpos,\ypos){\vector(\run,\rise){\arrowlength}}%
\or \put(\xpos,\ypos){\mvector(\run,\rise){\arrowlength}}%
\or \put(\xpos,\ypos){\evector(\run,\rise){\arrowlength}}%
\fi}\fi\fi\fi\fi}}

\newcount\numbdashes \newcount\lengthdash \newcount\increment

\def\howmanydashes{
\numbdashes=\arrowlength \lengthdash=40 \divide\numbdashes by
\lengthdash \lengthdash=\arrowlength \divide\lengthdash by
\numbdashes
\increment=\lengthdash \multiply\lengthdash by 3
\divide\lengthdash by 5 }

\def\putdashvector(#1)(#2,#3)#4#5{%
\ifnum#3=0 \putdashhvector(#1){#4}#5 \else \ifnum#2=0
\putdashvvector(#1){#4}#5\fi\fi}

\def\putdashhvector(#1,#2)#3#4{{%
\arrowlength=#3 \howmanydashes
\multiput(#1,#2)(\increment,0){\numbdashes}%
{\vrule height .4pt width \lengthdash\unitlength} \arrowtype=#4
\xpos=#1 \ifnum\arrowtype<0 \advance\arrowtype by 7 \fi
\ifcase\arrowtype \or \advance\xpos by 10
    \put(\xpos,#2){\vector(-1,0){\lengthdash}}
    \advance\xpos by 40
    \put(\xpos,#2){\vector(-1,0){\lengthdash}}
\or \advance \xpos by 10
    \put(\xpos,#2){\vector(-1,0){\lengthdash}}
    \advance\xpos by  \arrowlength
    \advance\xpos by  -50
    \put(\xpos,#2){\vector(-1,0){\lengthdash}}
\or \advance\xpos by 10
    \put(\xpos,#2){\vector(-1,0){\lengthdash}}
\or \advance\xpos by \arrowlength
    \advance\xpos by -\lengthdash
    \put(\xpos,#2){\vector(1,0){\lengthdash}}
\or {\advance\xpos by 10
    \put(\xpos,#2){\vector(1,0){\lengthdash}}}
    \advance\xpos by \arrowlength
    \advance\xpos by -\lengthdash
    \put(\xpos,#2){\vector(1,0){\lengthdash}}
\or \advance\xpos by \arrowlength
    \advance\xpos by -\lengthdash
    \put(\xpos,#2){\vector(1,0){\lengthdash}}
    \advance\xpos by -40
    \put(\xpos,#2){\vector(1,0){\lengthdash}}
   \fi
}}

\def\putdashvvector(#1,#2)#3#4{{%
\arrowlength=#3 \howmanydashes \ypos=#2 \advance\ypos by
-\arrowlength
\multiput(#1,#2)(0,\increment){\numbdashes}%
    {\vrule width .4pt height \lengthdash\unitlength}
\arrowtype=#4 \ypos=#2 \ifnum\arrowtype<0 \advance\arrowtype by 7
\fi \ifcase\arrowtype \or \advance\ypos by \arrowlength
\advance\ypos by -40
    \put(#1,\ypos){\vector(0,1){\lengthdash}}
    \advance\ypos by -40
    \put(#1,\ypos){\vector(0,1){\lengthdash}}
\or \advance\ypos by 10
    \put(#1,\ypos){\vector(0,1){\lengthdash}}
    \advance\ypos by \arrowlength \advance\ypos by -40
    \put(#1,\ypos){\vector(0,1){\lengthdash}}
\or \advance\ypos by \arrowlength \advance\ypos by -40
    \put(#1,\ypos){\vector(0,1){\lengthdash}}
\or \advance\ypos by 10
    \put(#1,\ypos){\vector(0,-1){\lengthdash}}
\or \advance\ypos by 10
    \put(#1,\ypos){\vector(0,-1){\lengthdash}}
    \advance\ypos by \arrowlength \advance\ypos by -40
    \put(#1,\ypos){\vector(0,-1){\lengthdash}}
\or \advance\ypos by 10
    \put(#1,\ypos){\vector(0,-1){\lengthdash}}
    \advance\ypos by 40
    \put(#1,\ypos){\vector(0,-1){\lengthdash}}
\fi }}

\def\puthmorphism(#1,#2)[#3`#4`#5]#6#7#8{{%
\xpos #1 \ypos #2 \width #6 \arrowlength #6 \arrowtype=#7
\putbox(\xpos,\ypos){#3\vphantom{#4}}%
{\advance \xpos by\arrowlength
\putbox(\xpos,\ypos){\vphantom{#3}#4}}%
\horsize{\tempcounta}{#3}%
\horsize{\tempcountb}{#4}%
\divide \tempcounta by2 \divide \tempcountb by2 \advance
\tempcounta by30 \advance \tempcountb by30 \advance \xpos
by\tempcounta \advance \arrowlength by-\tempcounta \advance
\arrowlength by-\tempcountb
\putvector(\xpos,\ypos)(1,0)\arrowlength\arrowtype \divide
\arrowlength by2 \advance \xpos by\arrowlength
\vertsize{\tempcounta}{#5}%
\divide\tempcounta by2 \advance \tempcounta by20
\if a#8 %
   \advance \ypos by\tempcounta
   \putbox(\xpos,\ypos){#5}%
\else
   \advance \ypos by-\tempcounta
   \putbox(\xpos,\ypos){#5}%
\fi}}

\def\putvmorphism(#1,#2)[#3`#4`#5]#6#7#8{{%
\xpos #1 \ypos #2 \arrowlength #6 \arrowtype #7
\settowidth{\xlen}{$#5$}%
\putbox(\xpos,\ypos){#3}%
{\advance \ypos by-\arrowlength
\putbox(\xpos,\ypos){#4}}%
{\advance\arrowlength by-140 \advance \ypos by-70 \ifdim\xlen>0pt
   \if m#8%
      \putsplitvector(\xpos,\ypos)\arrowlength\arrowtype
   \else
   \putvector(\xpos,\ypos)(0,-1)\arrowlength\arrowtype
   \fi
\else
   \putvector(\xpos,\ypos)(0,-1)\arrowlength\arrowtype
\fi}%
\ifdim\xlen>0pt
   \divide \arrowlength by2
   \advance\ypos by-\arrowlength
   \if l#8%
      \advance \xpos by-40
      \putrbox(\xpos,\ypos){#5}%
   \else\if r#8%
      \advance \xpos by40
      \putlbox(\xpos,\ypos){#5}%
   \else
      \putbox(\xpos,\ypos){#5}%
   \fi\fi
\fi }}

\def\putsquarep<#1>(#2)[#3;#4`#5`#6`#7]{{%
\setsqparms[#1]%
\setpos(#2)%
\settokens`#3`%
\puthmorphism(\xpos,\ypos)[\tokenc`\tokend`{#7}]{\width}{\arrowtyped}b%
\advance\ypos by \height
\puthmorphism(\xpos,\ypos)[\tokena`\tokenb`{#4}]{\width}{\arrowtypea}a%
\putvmorphism(\xpos,\ypos)[``{#5}]{\height}{\arrowtypeb}l%
\advance\xpos by \width
\putvmorphism(\xpos,\ypos)[``{#6}]{\height}{\arrowtypec}r%
}}

\def\putsquare{\@ifnextchar <{\putsquarep}{\putsquarep%
   <\arrowtypea`\arrowtypeb`\arrowtypec`\arrowtyped;\width`\height>}}
\def\square{\@ifnextchar< {\squarep}{\squarep
   <\arrowtypea`\arrowtypeb`\arrowtypec`\arrowtyped;\width`\height>}}
\def\squarep<#1>[#2`#3`#4`#5;#6`#7`#8`#9]{{
\setsqparms[#1]
\diagram
\putsquarep<\arrowtypea`\arrowtypeb`\arrowtypec`
\arrowtyped;\width`\height>
(0,0)[#2`#3`#4`{#5};#6`#7`#8`{#9}]
\enddiagram
}}                                                 
\def\putptrianglep<#1>(#2,#3)[#4`#5`#6;#7`#8`#9]{{%
\settriparms[#1]%
\xpos=#2 \ypos=#3 \advance\ypos by \height
\puthmorphism(\xpos,\ypos)[#4`#5`{#7}]{\height}{\arrowtypea}a%
\putvmorphism(\xpos,\ypos)[`#6`{#8}]{\height}{\arrowtypeb}l%
\advance\xpos by\height
\putmorphism(\xpos,\ypos)(-1,-1)[``{#9}]{\height}{\arrowtypec}r%
}}

\def\putptriangle{\@ifnextchar <{\putptrianglep}{\putptrianglep
   <\arrowtypea`\arrowtypeb`\arrowtypec;\height>}}
\def\ptriangle{\@ifnextchar <{\ptrianglep}{\ptrianglep
   <\arrowtypea`\arrowtypeb`\arrowtypec;\height>}}
\def\ptrianglep<#1>[#2`#3`#4;#5`#6`#7]{{
\settriparms[#1]
\diagram
\putptrianglep<\arrowtypea`\arrowtypeb`
\arrowtypec;\height>
(0,0)[#2`#3`#4;#5`#6`{#7}]
\enddiagram
}}                                            

\def\putqtrianglep<#1>(#2,#3)[#4`#5`#6;#7`#8`#9]{{%
\settriparms[#1]%
\xpos=#2 \ypos=#3 \advance\ypos by\height
\puthmorphism(\xpos,\ypos)[#4`#5`{#7}]{\height}{\arrowtypea}a%
\putmorphism(\xpos,\ypos)(1,-1)[``{#8}]{\height}{\arrowtypeb}l%
\advance\xpos by\height
\putvmorphism(\xpos,\ypos)[`#6`{#9}]{\height}{\arrowtypec}r%
}}

\def\putqtriangle{\@ifnextchar <{\putqtrianglep}{\putqtrianglep
   <\arrowtypea`\arrowtypeb`\arrowtypec;\height>}}
\def\qtriangle{\@ifnextchar <{\qtrianglep}{\qtrianglep
   <\arrowtypea`\arrowtypeb`\arrowtypec;\height>}}
\def\qtrianglep<#1>[#2`#3`#4;#5`#6`#7]{{
\settriparms[#1]
\width=\height                                
\diagram
\putqtrianglep<\arrowtypea`\arrowtypeb`
\arrowtypec;\height>
(0,0)[#2`#3`#4;#5`#6`{#7}]
\enddiagram
}}

\def\putdtrianglep<#1>(#2,#3)[#4`#5`#6;#7`#8`#9]{{%
\settriparms[#1]%
\xpos=#2 \ypos=#3
\puthmorphism(\xpos,\ypos)[#5`#6`{#9}]{\height}{\arrowtypec}b%
\advance\xpos by \height \advance\ypos by\height
\putmorphism(\xpos,\ypos)(-1,-1)[``{#7}]{\height}{\arrowtypea}l%
\putvmorphism(\xpos,\ypos)[#4``{#8}]{\height}{\arrowtypeb}r%
}}

\def\putdtriangle{\@ifnextchar <{\putdtrianglep}{\putdtrianglep
   <\arrowtypea`\arrowtypeb`\arrowtypec;\height>}}
\def\dtriangle{\@ifnextchar <{\dtrianglep}{\dtrianglep
   <\arrowtypea`\arrowtypeb`\arrowtypec;\height>}}
\def\dtrianglep<#1>[#2`#3`#4;#5`#6`#7]{{
\settriparms[#1]
\width=\height                                
\diagram
\putdtrianglep<\arrowtypea`\arrowtypeb`
\arrowtypec;\height>
(0,0)[#2`#3`#4;#5`#6`{#7}]
\enddiagram
}}

\def\putbtrianglep<#1>(#2,#3)[#4`#5`#6;#7`#8`#9]{{%
\settriparms[#1]%
\xpos=#2 \ypos=#3
\puthmorphism(\xpos,\ypos)[#5`#6`{#9}]{\height}{\arrowtypec}b%
\advance\ypos by\height
\putmorphism(\xpos,\ypos)(1,-1)[``{#8}]{\height}{\arrowtypeb}r%
\putvmorphism(\xpos,\ypos)[#4``{#7}]{\height}{\arrowtypea}l%
}}

\def\putbtriangle{\@ifnextchar <{\putbtrianglep}{\putbtrianglep
   <\arrowtypea`\arrowtypeb`\arrowtypec;\height>}}
\def\btriangle{\@ifnextchar <{\btrianglep}{\btrianglep
   <\arrowtypea`\arrowtypeb`\arrowtypec;\height>}}
\def\btrianglep<#1>[#2`#3`#4;#5`#6`#7]{{
\settriparms[#1]
\width=\height                               
\diagram
\putbtrianglep<\arrowtypea`\arrowtypeb`
\arrowtypec;\height>
(0,0)[#2`#3`#4;#5`#6`{#7}]
\enddiagram
}}

\def\putAtrianglep<#1>(#2,#3)[#4`#5`#6;#7`#8`#9]{{%
\settriparms[#1]%
\xpos=#2 \ypos=#3 {\multiply \height by2
\puthmorphism(\xpos,\ypos)[#5`#6`{#9}]{\height}{\arrowtypec}b}%
\advance\xpos by\height \advance\ypos by\height
\putmorphism(\xpos,\ypos)(-1,-1)[#4``{#7}]{\height}{\arrowtypea}l%
\putmorphism(\xpos,\ypos)(1,-1)[``{#8}]{\height}{\arrowtypeb}r%
}}

\def\putAtriangle{\@ifnextchar <{\putAtrianglep}{\putAtrianglep
   <\arrowtypea`\arrowtypeb`\arrowtypec;\height>}}
\def\Atriangle{\@ifnextchar <{\Atrianglep}{\Atrianglep
   <\arrowtypea`\arrowtypeb`\arrowtypec;\height>}}
\def\Atrianglep<#1>[#2`#3`#4;#5`#6`#7]{{
\settriparms[#1]
\width=\height                                     
\diagram
\putAtrianglep<\arrowtypea`\arrowtypeb`
\arrowtypec;\height>
(0,0)[#2`#3`#4;#5`#6`{#7}]
\enddiagram
}}

\def\putAtrianglepairp<#1>(#2)[#3;#4`#5`#6`#7`#8]{{%
\settripairparms[#1]%
\setpos(#2)%
\settokens`#3`%
\puthmorphism(\xpos,\ypos)[\tokenb`\tokenc`{#7}]{\height}{\arrowtyped}b%
\advance\xpos by\height
\puthmorphism(\xpos,\ypos)[\phantom{\tokenc}`\tokend`{#8}]%
{\height}{\arrowtypee}b%
\advance\ypos by\height
\putmorphism(\xpos,\ypos)(-1,-1)[\tokena``{#4}]{\height}{\arrowtypea}l%
\putvmorphism(\xpos,\ypos)[``{#5}]{\height}{\arrowtypeb}m%
\putmorphism(\xpos,\ypos)(1,-1)[``{#6}]{\height}{\arrowtypec}r%
}}

\def\putAtrianglepair{\@ifnextchar <{\putAtrianglepairp}{\putAtrianglepairp%
   <\arrowtypea`\arrowtypeb`\arrowtypec`\arrowtyped`\arrowtypee;\height>}}
\def\Atrianglepair{\@ifnextchar <{\Atrianglepairp}{\Atrianglepairp%
   <\arrowtypea`\arrowtypeb`\arrowtypec`\arrowtyped`\arrowtypee;\height>}}

\def\Atrianglepairp<#1>[#2;#3`#4`#5`#6`#7]{{
\settripairparms[#1]
\settokens`#2`
\width=\height                                
\diagram
\putAtrianglepairp                            
<\arrowtypea`\arrowtypeb`\arrowtypec`
\arrowtyped`\arrowtypee;\height>
(0,0)[{#2};#3`#4`#5`#6`{#7}]
\enddiagram
}}

\def\putVtrianglep<#1>(#2,#3)[#4`#5`#6;#7`#8`#9]{{%
\settriparms[#1]%
\xpos=#2 \ypos=#3 \advance\ypos by\height {\multiply\height by2
\puthmorphism(\xpos,\ypos)[#4`#5`{#7}]{\height}{\arrowtypea}a}%
\putmorphism(\xpos,\ypos)(1,-1)[`#6`{#8}]{\height}{\arrowtypeb}l%
\advance\xpos by\height \advance\xpos by\height
\putmorphism(\xpos,\ypos)(-1,-1)[``{#9}]{\height}{\arrowtypec}r%
}}

\def\putVtriangle{\@ifnextchar <{\putVtrianglep}{\putVtrianglep
   <\arrowtypea`\arrowtypeb`\arrowtypec;\height>}}
\def\Vtriangle{\@ifnextchar <{\Vtrianglep}{\Vtrianglep
   <\arrowtypea`\arrowtypeb`\arrowtypec;\height>}}
\def\Vtrianglep<#1>[#2`#3`#4;#5`#6`#7]{{
\settriparms[#1]
\width=\height                                 
\diagram
\putVtrianglep<\arrowtypea`\arrowtypeb`
\arrowtypec;\height>
(0,0)[#2`#3`#4;#5`#6`{#7}]
\enddiagram
}}

\def\putVtrianglepairp<#1>(#2)[#3;#4`#5`#6`#7`#8]{{
\settripairparms[#1]%
\setpos(#2)%
\settokens`#3`%
\advance\ypos by\height
\putmorphism(\xpos,\ypos)(1,-1)[`\tokend`{#6}]{\height}{\arrowtypec}l%
\puthmorphism(\xpos,\ypos)[\tokena`\tokenb`{#4}]{\height}{\arrowtypea}a%
\advance\xpos by\height
\puthmorphism(\xpos,\ypos)[\phantom{\tokenb}`\tokenc`{#5}]%
{\height}{\arrowtypeb}a%
\putvmorphism(\xpos,\ypos)[``{#7}]{\height}{\arrowtyped}m%
\advance\xpos by\height
\putmorphism(\xpos,\ypos)(-1,-1)[``{#8}]{\height}{\arrowtypee}r%
}}

\def\putVtrianglepair{\@ifnextchar <{\putVtrianglepairp}{\putVtrianglepairp%
    <\arrowtypea`\arrowtypeb`\arrowtypec`\arrowtyped`\arrowtypee;\height>}}
\def\Vtrianglepair{\@ifnextchar <{\Vtrianglepairp}{\Vtrianglepairp%
    <\arrowtypea`\arrowtypeb`\arrowtypec`\arrowtyped`\arrowtypee;\height>}}
\def\Vtrianglepairp<#1>[#2;#3`#4`#5`#6`#7]{{
\settripairparms[#1]
\settokens`#2`
\diagram
\putVtrianglepairp                             
<\arrowtypea`\arrowtypeb`\arrowtypec`
\arrowtyped`\arrowtypee;\height>
(0,0)[{#2};#3`#4`#5`#6`{#7}]
\enddiagram
}}

\def\putCtrianglep<#1>(#2,#3)[#4`#5`#6;#7`#8`#9]{{%
\settriparms[#1]%
\xpos=#2 \ypos=#3 \advance\ypos by\height
\putmorphism(\xpos,\ypos)(1,-1)[``{#9}]{\height}{\arrowtypec}l%
\advance\xpos by\height \advance\ypos by\height
\putmorphism(\xpos,\ypos)(-1,-1)[#4`#5`{#7}]{\height}{\arrowtypea}l%
{\multiply\height by 2
\putvmorphism(\xpos,\ypos)[`#6`{#8}]{\height}{\arrowtypeb}r}%
}}

\def\putCtriangle{\@ifnextchar <{\putCtrianglep}{\putCtrianglep
    <\arrowtypea`\arrowtypeb`\arrowtypec;\height>}}
\def\Ctriangle{\@ifnextchar <{\Ctrianglep}{\Ctrianglep
    <\arrowtypea`\arrowtypeb`\arrowtypec;\height>}}
\def\Ctrianglep<#1>[#2`#3`#4;#5`#6`#7]{{
\settriparms[#1]
\width=\height                               
\diagram
\putCtrianglep<\arrowtypea`\arrowtypeb`
\arrowtypec;\height>
(0,0)[#2`#3`#4;#5`#6`{#7}]
\enddiagram
}}                                           
\def\putDtrianglep<#1>(#2,#3)[#4`#5`#6;#7`#8`#9]{{%
\settriparms[#1]%
\xpos=#2 \ypos=#3 \advance\xpos by\height \advance\ypos by\height
\putmorphism(\xpos,\ypos)(-1,-1)[``{#9}]{\height}{\arrowtypec}r%
\advance\xpos by-\height \advance\ypos by\height
\putmorphism(\xpos,\ypos)(1,-1)[`#5`{#8}]{\height}{\arrowtypeb}r%
{\multiply\height by 2
\putvmorphism(\xpos,\ypos)[#4`#6`{#7}]{\height}{\arrowtypea}l}%
}}

\def\putDtriangle{\@ifnextchar <{\putDtrianglep}{\putDtrianglep
    <\arrowtypea`\arrowtypeb`\arrowtypec;\height>}}
\def\Dtriangle{\@ifnextchar <{\Dtrianglep}{\Dtrianglep
   <\arrowtypea`\arrowtypeb`\arrowtypec;\height>}}
\def\Dtrianglep<#1>[#2`#3`#4;#5`#6`#7]{{
\settriparms[#1]
\width=\height                              
\diagram
\putDtrianglep<\arrowtypea`\arrowtypeb`
\arrowtypec;\height>
(0,0)[#2`#3`#4;#5`#6`{#7}]
\enddiagram
}}                                          
\def\setrecparms[#1`#2]{\width=#1 \height=#2}%

\def\recursep<#1`#2>[#3;#4`#5`#6`#7`#8]{{\m@th
\width=#1 \height=#2 \settokens`#3`
\settowidth{\tempdimen}{$\tokena$} \ifdim\tempdimen=0pt
  \savebox{\tempboxa}{\hbox{$\tokenb$}}%
  \savebox{\tempboxb}{\hbox{$\tokend$}}%
  \savebox{\tempboxc}{\hbox{$#6$}}%
\else
  \savebox{\tempboxa}{\hbox{$\hbox{$\tokena$}\times\hbox{$\tokenb$}$}}%
  \savebox{\tempboxb}{\hbox{$\hbox{$\tokena$}\times\hbox{$\tokend$}$}}%
  \savebox{\tempboxc}{\hbox{$\hbox{$\tokena$}\times\hbox{$#6$}$}}%
\fi \ypos=\height \divide\ypos by 2 \xpos=\ypos \advance\xpos by
\width \bfig
\putCtrianglep<-1`1`1;\ypos>(0,0)[`\tokenc`;#5`#6`{#7}]%
\puthmorphism(\ypos,0)[\tokend`\usebox{\tempboxb}`{#8}]{\width}{-1}b%
\puthmorphism(\ypos,\height)[\tokenb`\usebox{\tempboxa}`{#4}]{\width}{-1}a%
\advance\ypos by \width
\putvmorphism(\ypos,\height)[``\usebox{\tempboxc}]{\height}1r%
\efig }}

\def\recurse{\@ifnextchar <{\recursep}{\recursep<\width`\height>}}

\def\puttwohmorphisms(#1,#2)[#3`#4;#5`#6]#7#8#9{{%
\puthmorphism(#1,#2)[#3`#4`]{#7}0a \ypos=#2 \advance\ypos by 20
\puthmorphism(#1,\ypos)[\phantom{#3}`\phantom{#4}`#5]{#7}{#8}a
\advance\ypos by -40
\puthmorphism(#1,\ypos)[\phantom{#3}`\phantom{#4}`#6]{#7}{#9}b }}

\def\puttwovmorphisms(#1,#2)[#3`#4;#5`#6]#7#8#9{{%
\putvmorphism(#1,#2)[#3`#4`]{#7}0a \xpos=#1 \advance\xpos by -20
\putvmorphism(\xpos,#2)[\phantom{#3}`\phantom{#4}`#5]{#7}{#8}l
\advance\xpos by 40
\putvmorphism(\xpos,#2)[\phantom{#3}`\phantom{#4}`#6]{#7}{#9}r }}

\def\puthcoequalizer(#1)[#2`#3`#4;#5`#6`#7]#8#9{{%
\setpos(#1)%
\puttwohmorphisms(\xpos,\ypos)[#2`#3;#5`#6]{#8}11%
\advance\xpos by #8
\puthmorphism(\xpos,\ypos)[\phantom{#3}`#4`#7]{#8}1{#9} }}

\def\putvcoequalizer(#1)[#2`#3`#4;#5`#6`#7]#8#9{{%
\setpos(#1)%
\puttwovmorphisms(\xpos,\ypos)[#2`#3;#5`#6]{#8}11%
\advance\ypos by -#8
\putvmorphism(\xpos,\ypos)[\phantom{#3}`#4`#7]{#8}1{#9} }}

\def\putthreehmorphisms(#1)[#2`#3;#4`#5`#6]#7(#8)#9{{%
\setpos(#1) \settypes(#8)
\if a#9 %
     \vertsize{\tempcounta}{#5}%
     \vertsize{\tempcountb}{#6}%
     \ifnum \tempcounta<\tempcountb \tempcounta=\tempcountb \fi
\else
     \vertsize{\tempcounta}{#4}%
     \vertsize{\tempcountb}{#5}%
     \ifnum \tempcounta<\tempcountb \tempcounta=\tempcountb \fi
\fi \advance \tempcounta by 60
\puthmorphism(\xpos,\ypos)[#2`#3`#5]{#7}{\arrowtypeb}{#9}
\advance\ypos by \tempcounta
\puthmorphism(\xpos,\ypos)[\phantom{#2}`\phantom{#3}`#4]{#7}{\arrowtypea}{#9}
\advance\ypos by -\tempcounta \advance\ypos by -\tempcounta
\puthmorphism(\xpos,\ypos)[\phantom{#2}`\phantom{#3}`#6]{#7}{\arrowtypec}{#9}
}}

\def\setarrowtoks[#1`#2`#3`#4`#5`#6]{%
\def\toka{#1}
\def\tokb{#2}
\def\tokc{#3}
\def\tokd{#4}
\def\toke{#5}
\def\tokf{#6}
}
\def\hex{\@ifnextchar <{\hexp}{\hexp<1000`400>}}
\def\hexp<#1`#2>[#3`#4`#5`#6`#7`#8;#9]{%
\setarrowtoks[#9] \yext=#2 \advance \yext by #2 \xext=#1
\advance\xext by \yext \bfig
\putCtriangle<-1`0`1;#2>(0,0)[`#5`;\tokb``\tokd] \xext=#1
\yext=#2 \advance \yext by #2
\putsquare<1`0`0`1;\xext`\yext>(#2,0)[#3`#4`#7`#8;\toka```\tokf]
\advance \xext by #2
\putDtriangle<0`1`-1;#2>(\xext,0)[`#6`;`\tokc`\toke] \efig }

\makeatother

\input{tcilatex}
\begin{document}

\title{Undergraduate Lecture Notes in De Rham--Hodge Theory}
\author{Vladimir G. Ivancevic\thanks{%
Land Operations Division, Defence Science \& Technology Organisation, P.O.
Box 1500, Edinburgh SA 5111, Australia
(Vladimir.Ivancevic@dsto.defence.gov.au)} \and Tijana T. Ivancevic\thanks{
Tesla Science Evolution Institute, Adelaide, Australia (Tijana.Ivancevic@alumni.adelaide.edu.au)}}
\date{}
\maketitle

\begin{abstract}
These lecture notes in the De Rham--Hodge theory are designed for a
1--semester undergraduate course (in mathematics, physics, engineering,
chemistry or biology). This landmark theory of the 20th Century mathematics
gives a rigorous foundation to modern field and gauge theories in physics,
engineering and physiology. The only necessary background for comprehensive
reading of these notes is Green's theorem from multivariable calculus.
\end{abstract}

\tableofcontents

\section{Exterior Geometrical Mmachinery}

To grasp the essence of Hodge--De Rham theory, we need first to familiarize
ourselves with exterior differential forms and Stokes' theorem.

\subsection{From Green's to Stokes' theorem}

Recall that \textit{Green's theorem} in the region $C\ $in $(x,y)-$plane $%
\mathbb{R}^{2}$\ connects a line integral $\oint_{\partial C}$ (over the
boundary $\partial C$ of $C)$ with a double integral $\iint_{C}$ over $C$
(see e.g., \cite{MarsVec})
\begin{equation*}
\oint_{\partial C}Pdx+Qdy=\iint_{C}\left( \frac{\partial Q}{\partial x}-%
\frac{\partial P}{\partial y}\right) dxdy.
\end{equation*}%
In other words, if we define two differential forms (integrands of $%
\oint_{\partial C}$ and $\iint_{C}$) as
\begin{eqnarray*}
\text{1--form} &:&\mathbf{A}=Pdx+Qdy,\qquad \text{and} \\
\text{2--form} &:&\mathbf{dA}=\left( \frac{\partial Q}{\partial x}-\frac{%
\partial P}{\partial y}\right) dxdy,
\end{eqnarray*}%
(where $\mathbf{d}$ denotes the \textit{exterior derivative} that makes a $%
(p+1)-$form out of a $p-$form, see next subsection), then we can rewrite
Green's theorem as \textit{Stokes' theorem:}
\begin{equation*}
\int_{\partial C}\mathbf{A}=\int_{C}\mathbf{dA}.
\end{equation*}%
The integration domain $C$ is in topology called a \textit{chain}, and $%
\partial C$ is a 1D boundary of a 2D chain $C$. In general, \emph{the
boundary of a boundary is zero} (see \cite{MTW,CW}), that is, $\partial
(\partial C)=0$, or formally $\partial ^{2}=0$.

\subsection{Exterior derivative}

\begin{figure}[tbp]
\centerline{\includegraphics[height=13cm]{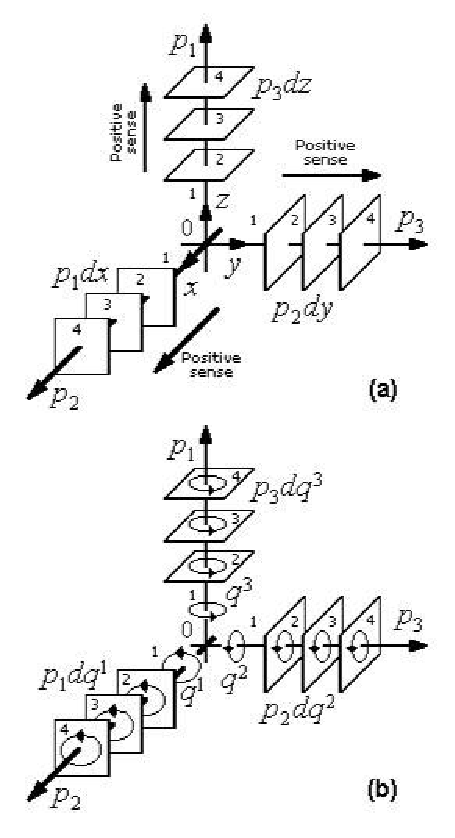}}
\caption{Basis vectors and one-forms in Euclidean $\mathbb{R}^{3}-$space:
(a) Translational case; and (b) Rotational case \protect\cite{VladSIAM}. For
the same geometry in $\mathbb{R}^{3}$, see \protect\cite{MTW}.}
\label{Basis}
\end{figure}

\begin{figure}[tbp]
\centerline{\includegraphics[height=13cm]{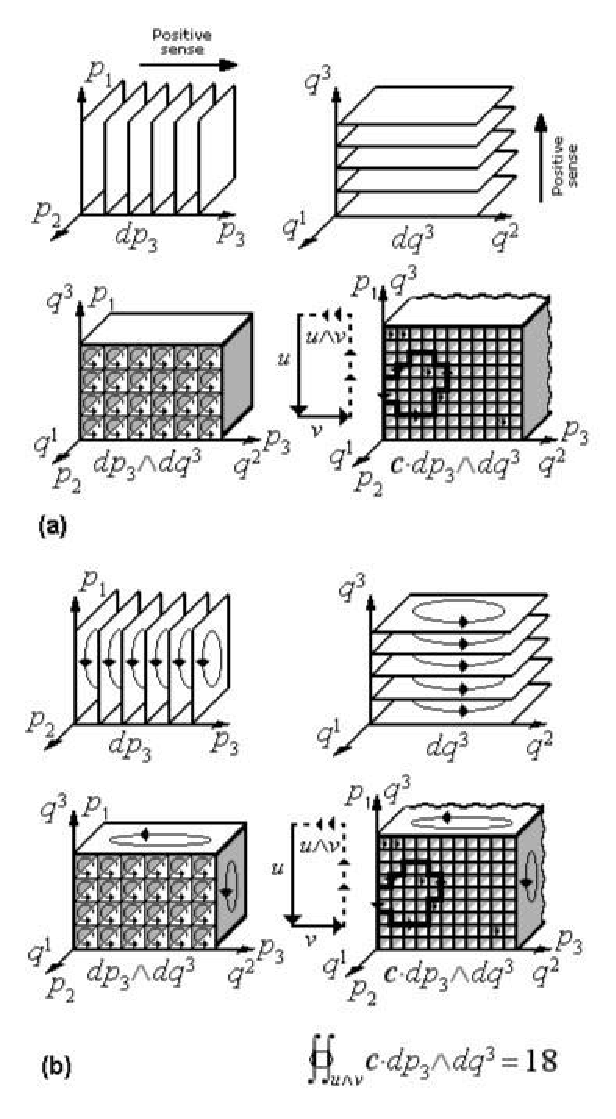}}
\caption{Fundamental two--form and its flux in $\mathbb{R}^3$: (a)
Translational case; (b) Rotational case. In both cases the flux through the
plane $u\wedge v$ is defined as $\protect\int\protect\int_{u\wedge v}
c\,dp_idq^i$ and measured by the number of tubes crossed by the circulation
oriented by $u\wedge v$ \protect\cite{VladSIAM}. For the same geometry in $%
\mathbb{R}^3$, see \protect\cite{MTW}.}
\label{2form}
\end{figure}

The exterior derivative $\mathbf{d}$ is a generalization of ordinary vector
differential operators (\textsl{grad}, \textsl{div} and \textsl{curl} see
\cite{De Rham,Flanders}) that transforms $p-$forms $\mathbf{\omega }$ into $%
(p+1)-$forms $\mathbf{d\omega }$ (see next subsection), with the main
property: $\mathbf{dd}=\mathbf{d}^{2}=0$, so that in $\mathbb{R}^{3}$ we
have (see Figures \ref{Basis} and \ref{2form})

\begin{itemize}
\item any scalar function $f=f(x,y,z)$ is a 0--form;

\item the gradient $\mathbf{d}f=\mathbf{\omega }$ of any smooth function $f$
is a 1--form
\begin{equation*}
\mathbf{\omega =d}f=\frac{\partial f}{\partial x}dx+\frac{\partial f}{%
\partial y}dy+\frac{\partial f}{\partial z}dz;
\end{equation*}

\item the curl $\mathbf{\alpha =d\omega}$ of any smooth 1--form $\mathbf{%
\omega }$ is a 2--form
\begin{eqnarray*}
&&\mathbf{\alpha =d\omega} =\left( \frac{\partial R}{\partial y}-\frac{%
\partial Q}{\partial z}\right) dydz+\left( \frac{\partial P}{\partial z}-%
\frac{\partial R}{\partial x}\right) dzdx+\left( \frac{\partial Q}{\partial x%
}-\frac{\partial P}{\partial y}\right) dxdy; \\
&&\text{if }\mathbf{\omega =d}f\Rightarrow \mathbf{\alpha =dd}f=0.
\end{eqnarray*}

\item the divergence $\mathbf{\beta =d\alpha}$ of any smooth 2--form $%
\mathbf{\alpha } $ is a 3--form
\begin{equation*}
\mathbf{\beta =d\alpha} =\left( \frac{\partial A}{\partial x}+\frac{\partial
B}{\partial y}+\frac{\partial C}{\partial z}\right) dxdydz;\qquad\text{if }%
\mathbf{\alpha =d\omega}\Rightarrow \mathbf{\beta =dd\omega}=0.
\end{equation*}
\end{itemize}

In general, for any two smooth functions $f=f(x,y,z)$ and $g=g(x,y,z)$, the
exterior derivative $\mathbf{d}$\ obeys the \emph{Leibniz rule} \cite%
{GaneshSprBig,GaneshADG}:
\begin{equation*}
\mathbf{d}(fg)=g\,\mathbf{d}f+f\,\mathbf{d}g,
\end{equation*}%
and the \textit{chain rule}:
\begin{equation*}
\mathbf{d}\left( g(f)\right) =g^{\prime }(f)\,\mathbf{d}f.
\end{equation*}

\subsection{Exterior forms}

In general, given a so--called 4D \textit{coframe}, that is a set of
coordinate differentials $\{dx^{i}\}\in \mathbb{R}^{4}$, we can define the
space of all $p-$forms, denoted $\Omega ^{p}(\mathbb{R}^{4})$, using the
exterior derivative $\mathbf{d}:\Omega ^{p}(\mathbb{R}^{4})\rightarrow
\Omega ^{p+1}(\mathbb{R}^{4})$ and Einstein's summation convention over
repeated indices (e.g., $A_{i}\,dx^{i}=\sum_{i=0}^{3}A_{i}\,dx^{i}$), we
have:

\begin{description}
\item[1--form] -- a generalization of the Green's 1--form $Pdx+Qdy$,
\begin{equation*}
\mathbf{A}=A_{i}\,dx^{i}\in \Omega ^{1}(\mathbb{R}^{4}).
\end{equation*}%
For example, in 4D electrodynamics, $\mathbf{A}$ represents electromagnetic
(co)vector potential.

\item[2--form] -- generalizing the Green's 2--form $(\partial _{x}Q-\partial
_{y}P)\,dxdy$ ~$($with $\partial _{j}=\partial /\partial x^{j}$),
\begin{eqnarray*}
\mathbf{B} &=&\mathbf{dA}\in \Omega ^{2}(\mathbb{R}^{4}),\qquad \text{with
components} \\
\mathbf{B} &=&\frac{1}{2}B_{ij}\,dx^{i}\wedge dx^{j},\qquad \text{or} \\
\mathbf{B} &=&\partial _{j}A_{i}\,dx^{j}\wedge dx^{i},\qquad \text{so that}
\\
B_{ij} &=&-2\partial _{j}A_{i}=\partial _{i}A_{j}-\partial _{j}A_{i}=-B_{ji}.
\end{eqnarray*}%
where $\wedge $ is the anticommutative exterior (or, `wedge') product of two
differential forms; given a $p-$form $\mathbf{\alpha }\in \Omega ^{p}(%
\mathbb{R}^{4})$ and a $q-$form $\mathbf{\beta }\in \Omega ^{q}(\mathbb{R}%
^{4}),$ their exterior product is a $(p+q)-$form $\mathbf{\alpha }\wedge
\mathbf{\beta }\in \Omega ^{p+q}(\mathbb{R}^{4})$; e.g., if we have two
1--forms $\mathbf{a}=a_{i}dx^{i}$, and $\mathbf{b}=b_{i}dx^{i}$, their wedge
product $\mathbf{a}\wedge \mathbf{b}$ is a 2--form $\mathbf{\alpha }$ given
by
\begin{equation*}
\mathbf{\alpha }=\mathbf{a}\wedge \mathbf{b}=a_{i}b_{j}\,dx^{i}\wedge
dx^{j}=-a_{i}b_{j}\,dx^{j}\wedge dx^{i}=-\mathbf{b}\wedge \mathbf{a}.
\end{equation*}%
The exterior product $\wedge $ is related to the exterior derivative $%
\mathbf{d}=\partial_i dx^i$, by
\begin{equation*}
\mathbf{d}(\mathbf{\alpha }\wedge \mathbf{\beta })=\mathbf{d\alpha }\wedge
\mathbf{\beta }+(-1)^{p}\mathbf{\alpha }\wedge \mathbf{d\beta }.\text{ \ }
\end{equation*}

\item[3--form]
\begin{eqnarray*}
\mathbf{C} &=&\mathbf{dB~(}=\mathbf{ddA}\equiv 0)\in \Omega ^{3}(\mathbb{R}%
^{4}),\qquad \text{with components} \\
\mathbf{C} &=&\frac{1}{3!}C_{ijk}\,dx^{i}\wedge dx^{j}\wedge dx^{k},\qquad
\text{or} \\
\mathbf{C} &=&\partial _{k}B_{[ij]}\,dx^{k}\wedge dx^{i}\wedge dx^{j},\qquad
\text{so that} \\
C_{ijk} &=&-6\partial _{k}B_{[ij]},\qquad \text{where }B_{[ij]}\text{ is the
skew--symmetric part of }B_{ij}.
\end{eqnarray*}%
For example, in the 4D electrodynamics, $\mathbf{B}$ represents the field
2--form \emph{Faraday}, or the Li\'{e}nard--Wiechert 2--form (in the next
section we will use the standard symbol $\mathbf{F}$ instead of $\mathbf{B}$%
) satisfying the sourceless magnetic Maxwell's equation,
\begin{equation*}
\text{Bianchi identity}:\quad \mathbf{dB}=0,\quad \text{in components\quad }%
\partial _{k}B_{[ij]}=0.
\end{equation*}

\item[4--form]
\begin{eqnarray*}
\mathbf{D} &=&\mathbf{dC~(}=\mathbf{ddB}\equiv 0)\in \Omega ^{4}(\mathbb{R}%
^{4}),\qquad \text{with components} \\
\mathbf{D} &=&\partial _{l}C_{[ijk]}\,dx^{l}\wedge dx^{i}\wedge dx^{j}\wedge
dx^{k},\qquad \text{or} \\
\mathbf{D} &=&\frac{1}{4!}D_{ijkl}\,dx^{i}\wedge dx^{j}\wedge dx^{k}\wedge
dx^{l},\qquad \text{so that} \\
D_{ijkl} &=&-24\partial _{l}C_{[ijk]}.
\end{eqnarray*}
\end{description}

\subsection{Stokes theorem}

Generalization of the Green's theorem in the plane (and all other integral
theorems from vector calculus) is the Stokes theorem for the $p-$form $%
\mathbf{\omega }$, in an oriented $n$D domain $C$ (which is a $p-$chain with
a $(p-1)-$boundary $\partial C$, see next section)
\begin{equation*}
\int_{\partial C}\mathbf{\omega }=\int_{C}\mathbf{d\omega }.
\end{equation*}

For example, in the 4D Euclidean space $\mathbb{R}^{4}$ we have the
following three particular cases of the Stokes theorem, related to the
subspaces $C$ of $\mathbb{R}^{4}$:\newline
The 2D Stokes theorem:
\begin{equation*}
\int_{\partial C^{2}}\mathbf{A}=\int_{C^{2}}\mathbf{B}.
\end{equation*}
The 3D Stokes theorem:
\begin{equation*}
\int_{\partial C^{3}}\mathbf{B}=\int_{C^{3}}\mathbf{C}.
\end{equation*}
The 4D Stokes theorem:
\begin{equation*}
\int_{\partial C^{4}}\mathbf{C}=\int_{C^{4}}\mathbf{D}.
\end{equation*}

\section{De Rham--Hodge Theory Basics}

Now that we are familiar with differential forms and Stokes' theorem, we can
introduce Hodge--De Rham theory.

\subsection{Exact and closed forms and chains}

Notation change: we drop boldface letters from now on. In general, a $p-$%
form $\beta $ is called \textit{closed} if its exterior derivative $%
d=\partial_i dx^i$ is equal to zero,
\begin{equation*}
d\beta=0.
\end{equation*}
From this condition one can see that the closed form (the \textit{kernel} of
the exterior derivative operator $d$) is conserved quantity. Therefore,
closed $p-$forms possess certain invariant properties, physically
corresponding to the conservation laws (see e.g., \cite{Abraham}).

Also, a $p-$form $\beta$ that is an exterior derivative of some $(p-1)-$form
$\alpha$,
\begin{equation*}
\beta=d\alpha,
\end{equation*}
is called \textit{exact} (the \textit{image} of the exterior derivative
operator $d$). By Poincar\'e lemma, exact forms prove to be closed
automatically,
\begin{equation*}
d\beta=d(d\alpha)=0.
\end{equation*}

Since $d^{2}=0$, \emph{every exact form is closed.} The converse is only
partially true, by Poincar\'{e} lemma: every closed form is \textit{locally
exact}.

\begin{figure}[h]
\centerline{\includegraphics[width=8cm]{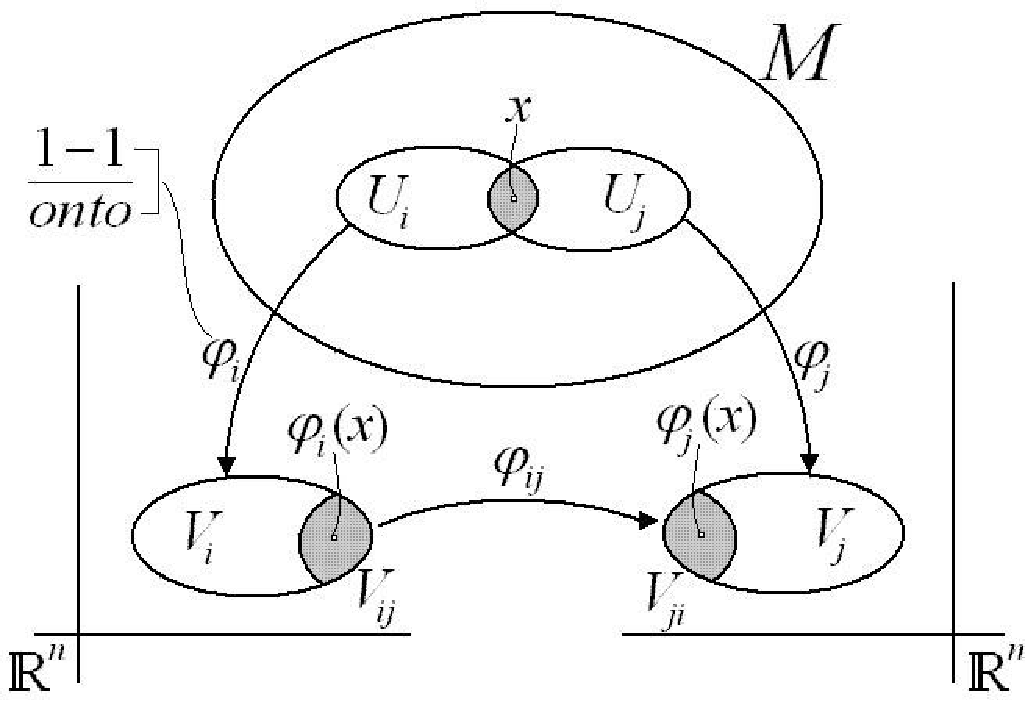}}
\caption{Depicting the manifold concept (see footnote 1).}
\label{Manifold1}
\end{figure}

Technically, this means that given a closed $p-$form $\alpha \in \Omega
^{p}(U)$, defined on an open set $U$ of a smooth manifold $M$ (see Figure %
\ref{Manifold1}), \footnote{%
Smooth manifold is a curved $n$D space which is locally equivalent to $%
\mathbb{R}^n$. To sketch it formal definition, consider a set $M$ (see
Figure \ref{Manifold1}) which is a \emph{candidate} for a manifold. Any
point $x\in M$ has its \textit{Euclidean chart}, given by a 1--1 and \emph{%
onto} map $\varphi _{i}:M\rightarrow \mathbb{R}^{n}$, with its \textit{%
Euclidean image} $V_{i}=\varphi_{i}(U_{i})$. Formally, a chart $\varphi_{i}$
is defined by
\begin{equation*}
\varphi _{i}:M\supset U_{i}\ni x\mapsto \varphi _{i}(x)\in V_{i}\subset
\mathbb{R}^{n},
\end{equation*}
where $U_{i}\subset M$ and $V_{i}\subset \mathbb{R}^{n}$ are open sets.
\par
Any point $x\in M$ can have several different charts (see Figure \ref%
{Manifold1}). Consider a case of two charts, $\varphi _{i},\varphi
_{j}:M\rightarrow \mathbb{R}^{n}$, having in their images two open sets, $%
V_{ij}=\varphi _{i}(U_{i}\cap U_{j})$ and $V_{ji}=\varphi _{j}(U_{i}\cap
U_{j})$. Then we have \textit{transition functions} $\varphi _{ij}$ between
them,
\begin{equation*}
\varphi _{ij}=\varphi _{j}\circ \varphi _{i}^{-1}:V_{ij}\rightarrow
V_{ji},\qquad \text{locally given by\qquad }\varphi _{ij}(x)=\varphi
_{j}(\varphi _{i}^{-1}(x)).
\end{equation*}
If transition functions $\varphi _{ij}$ exist, then we say that two charts, $%
\varphi _{i}$ and $\varphi _{j}$ are \emph{compatible}. Transition functions
represent a general (nonlinear) \emph{transformations of coordinates}, which
are the core of classical \emph{tensor calculus}.
\par
A set of compatible charts $\varphi _{i}:M\rightarrow \mathbb{R}^{n},$ such
that each point $x\in M$ has its Euclidean image in at least one chart, is
called an \textit{atlas}. Two atlases are \emph{equivalent} iff all their
charts are compatible (i.e., transition functions exist between them), so
their union is also an atlas. A \textit{manifold structure} is a class of
equivalent atlases.
\par
Finally, as charts $\varphi _{i}:M\rightarrow \mathbb{R}^{n}$ were supposed
to be 1-1 and onto maps, they can be either \emph{homeomorphism}\emph{s}, in
which case we have a \emph{topological} ($C^0$) manifold, or \emph{%
diffeomorphism}\emph{s}, in which case we have a \emph{smooth} ($C^{k}$)
manifold.} any point $m\in U$ has a neighborhood on which there exists a $%
(p-1)-$form $\beta \in \Omega ^{p-1}(U)$ such that $d\beta =\alpha |_{U}.$
In particular, there is a Poincar\'{e} lemma for contractible manifolds: Any
closed form on a smoothly contractible manifold is exact.

The Poincar\'{e} lemma is a generalization and unification of two
well--known facts in vector calculus:\newline
(i) If $\limfunc{curl}F=0$, then locally $F=\limfunc{grad}f$; ~ and ~ (ii)
If $\limfunc{div}F=0$, then locally $F=\limfunc{curl}G$.

A \textit{cycle} is a $p-$chain, (or, an oriented $p-$domain) $C\in \mathcal{%
C}_{p}(M)$ such that $\partial C=0 $. A \textit{boundary} is a chain $C$
such that $C=\partial B,$ for any other chain $B\in \mathcal{C}_{p}(M)$.
Similarly, a \textit{cocycle} (i.e., a \textit{closed form}) is a cochain $%
\omega $ such that $d\omega =0$. A \textit{coboundary} (i.e., an \textit{%
exact form}) is a cochain $\omega $ such that $\omega =d\theta ,$ for any
other cochain $\theta $. All exact forms are closed ($\omega =d\theta
\Rightarrow d\omega =0$) and all boundaries are cycles ($C=\partial
B\Rightarrow \partial C=0$). Converse is true only for smooth contractible
manifolds, by Poincar\'{e} lemma.

\subsection{De Rham duality of forms and chains}

Integration on a smooth manifold $M$ should be thought of as a nondegenerate
bilinear pairing $\left\langle , \right\rangle$ between $p-$forms and $p-$%
chains (spanning a finite domain on $M$). Duality of $p-$forms and $p-$%
chains on $M$ is based on the De Rham's `period', defined as \cite{{De
Rham,Bruhat}}
\begin{equation*}
\text{Period} :=\int_{C}\omega :=\left\langle C,\omega \right\rangle ,
\end{equation*}%
where $C$ is a cycle, $\omega $ is a cocycle, while $\left\langle C,\omega
\right\rangle =\omega (C)$ is their inner product $\left\langle C,\omega
\right\rangle :\Omega ^{p}(M)\times \mathcal{C}_{p}(M)\rightarrow \mathbb{R}$%
. From the Poincar\'{e} lemma, a closed $p-$form $\omega $ is exact iff $%
\left\langle C,\omega \right\rangle =0$.

The fundamental topological duality is based on the Stokes theorem,
\begin{equation*}
\int_{\partial C}\omega =\int_{C}d\omega \qquad \mathrm{or}\qquad
\left\langle \partial C,\omega \right\rangle =\left\langle C,d\omega
\right\rangle ,
\end{equation*}%
where $\partial C$ is the boundary of the $p-$chain $C$ oriented coherently
with $C$ on $M$. While the \textit{boundary operator} $\partial $ is a
global operator, the coboundary operator $d$ is local, and thus more
suitable for applications. The main property of the exterior differential,
\begin{equation*}
d\circ d\equiv d^{2}=0\quad \Longrightarrow \quad \partial \circ \partial
\equiv \partial ^{2}=0,\qquad (\text{and converse}),
\end{equation*}%
\noindent can be easily proved using the Stokes' theorem (and the above
`period notation') as
\begin{equation*}
0=\left\langle \partial ^{2}C,\omega \right\rangle =\left\langle \partial
C,d\omega \right\rangle =\left\langle C,d^{2}\omega \right\rangle =0.
\end{equation*}

\subsection{De Rham cochain and chain complexes}

In the Euclidean 3D space $\mathbb{R}^{3}$ we have the following De Rham
\emph{cochain complex}
\begin{equation*}
0\rightarrow \Omega ^{0}(\mathbb{R}^{3})\underset{\mathrm{grad}}{\overset{d}{%
\longrightarrow }}\Omega ^{1}(\mathbb{R}^{3})\underset{\mathrm{curl}}{%
\overset{d}{\longrightarrow }}\Omega ^{2}(\mathbb{R}^{3})\underset{\mathrm{%
div}}{\overset{d}{\longrightarrow }}\Omega ^{3}(\mathbb{R}^{3})\rightarrow 0.
\end{equation*}
Using the \textit{closure property} for the exterior differential in $%
\mathbb{R}^{3},~d\circ d\equiv d^{2}=0$, we get the standard identities from
vector calculus
\begin{equation*}
\limfunc{curl}\cdot \limfunc{grad}=0\text{ \ \ \ \ \ \ and \ \ \ \ \ \ }%
\limfunc{div}\cdot \limfunc{curl}=0.
\end{equation*}

As a duality, in $\mathbb{R}^{3}$ we have the following \emph{chain complex}
\begin{equation*}
0\leftarrow \mathcal{C}_{0}(\mathbb{R}^{3}){\overset{\partial }{%
\longleftarrow }}\mathcal{C}_{1}(\mathbb{R}^{3}){\overset{\partial }{%
\longleftarrow }}\mathcal{C}_{2}(\mathbb{R}^{3}){\overset{\partial }{%
\longleftarrow }}\mathcal{C}_{3}(\mathbb{R}^{3})\leftarrow 0,
\end{equation*}
(with the closure property $\partial\circ \partial\equiv \partial^{2}=0$)
which implies the following three boundaries:
\begin{equation*}
C_{1}\overset{\partial }{\mapsto }C_{0}=\partial (C_{1}),\qquad C_{2}\overset%
{\partial }{\mapsto }C_{1}=\partial (C_{2}),\qquad C_{3}\overset{\partial }{%
\mapsto }C_{2}=\partial (C_{3}),
\end{equation*}
where $C_0\in\mathcal{C}_0$ is a 0--boundary (or, a point), $C_1\in\mathcal{C%
}_1$ is a 1--boundary (or, a line), $C_2\in\mathcal{C}_2$ is a 2--boundary
(or, a surface), and $C_3\in\mathcal{C}_3$ is a 3--boundary (or, a
hypersurface). Similarly, the De Rham complex implies the following three
coboundaries:
\begin{equation*}
C^{0}\overset{d }{\mapsto }C^{1}=d (C^{0}),\qquad C^{1}\overset{d }{\mapsto }%
C^{2}=d (C^{1}),\qquad C^{2}\overset{d }{\mapsto }C^{3}=d (C^{2}),
\end{equation*}
where $C^0\in\Omega^0$ is 0--form (or, a function), $C^1\in\Omega^1$ is a
1--form, $C^2\in\Omega^2$ is a 2--form, and $C^3\in\Omega^3$ is a 3--form.

In general, on a smooth $n$D manifold $M$ we have the following De Rham
cochain complex \cite{De Rham}
\begin{equation*}
0\rightarrow \Omega ^{0}(M)\overset{d}{\longrightarrow }\Omega ^{1}(M)%
\overset{d}{\longrightarrow }\Omega ^{2}(M)\overset{d}{\longrightarrow }%
\Omega ^{3}(M)\overset{d}{\longrightarrow }\cdot \cdot \cdot \overset{d}{%
\longrightarrow }\Omega ^{n}(M)\rightarrow 0,
\end{equation*}%
satisfying the closure property on $M,~d\circ d\equiv d^{2}=0$.
\begin{figure}[tbh]
\centerline{\includegraphics[width=12cm]{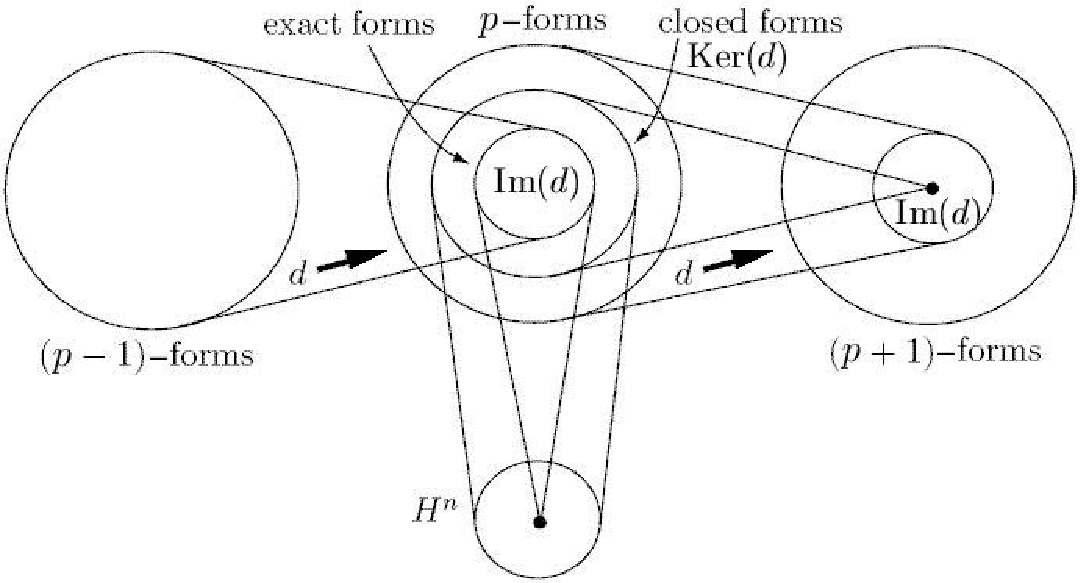}}
\caption{A small portion of the De Rham cochain complex, showing a
homomorphism of cohomology groups.}
\label{Cohomology}
\end{figure}

\subsection{De Rham cohomology vs. chain homology}

Briefly, the De Rham cohomology is the (functional) space of closed
differential $p-$forms modulo exact ones on a smooth manifold.

More precisely, the subspace of all closed $p-$forms (cocycles) on a smooth
manifold $M$ is the kernel $\limfunc{Ker}(d)$ of the De Rham $d-$%
homomorphism (see Figure \ref{Cohomology}), denoted by $Z^{p}(M)\subset
\Omega ^{p}(M)$, and the sub-subspace of all exact $p-$forms (coboundaries)
on $M$ is the image $\limfunc{Im}(d)$ of the De Rham homomorphism denoted by
$B^{p}(M)\subset Z^{p}(M)$. The \textit{quotient space}
\begin{equation}
H^{p}_{DR}(M):=\frac{Z^{p}(M)}{B^{p}{M}}=\frac{\limfunc{Ker}\left( d:\Omega
^{p}(M)\rightarrow \Omega ^{p+1}(M)\right) }{\limfunc{Im}\left( d:\Omega
^{p-1}(M)\rightarrow \Omega ^{p}(M)\right) },  \label{DR}
\end{equation}%
is called the $p$th De Rham \textit{cohomology group} of a manifold $M$. It
is a topological invariant of a manifold. Two $p-$cocycles $\alpha $,$\beta
\in \Omega ^{p}(M)$ are \emph{cohomologous}, or belong to the same \textit{%
cohomology class} $[\alpha ]\in H^{p}(M)$, if they differ by a $(p-1)-$%
coboundary $\alpha -\beta =d\theta \in \Omega ^{p-1}(M)$. The dimension $%
b_p=\dim H^{p}(M)$ of the De Rham cohomology group $H^{p}_{DR}(M)$ of the
manifold $M$ is called the Betti number $b_p$.

Similarly, the subspace of all $p-$cycles on a smooth manifold $M$ is the
kernel $\limfunc{Ker}(\partial )$ of the $\partial -$homomorphism, denoted
by $Z_{p}(M)\subset \mathcal{C}_{p}(M)$, and the sub-subspace of all $p-$%
boundaries on $M$ is the image $\limfunc{Im}(\partial )$ of the $\partial -$%
homomorphism, denoted by $B_{p}(M)\subset \mathcal{C}_{p}(M)$. Two $p-$%
cycles $C_{1}$,$C_{2}\in \mathcal{C}_{p}$ are \emph{homologous}, if they
differ by a $(p-1)-$boundary $C_{1}-C_{2}=\partial B\in \mathcal{C}_{p-1}(M)$%
. Then $C_{1}$ and $C_{2}$ belong to the same \textit{homology class} $%
[C]\in H_{p}(M)$, where $H_{p}(M)$ is the homology group of the manifold $M$%
, defined as
\begin{equation*}
H_{p}(M):=\frac{Z_{p}(M)}{B_{p}(M)}=\frac{\limfunc{Ker}(\partial :\mathcal{C}%
_{p}(M)\rightarrow \mathcal{C}_{p-1}(M))}{\limfunc{Im}(\partial :\mathcal{C}%
_{p+1}(M)\rightarrow \mathcal{C}_{p}(M))},
\end{equation*}%
where $Z_{p}$ is the vector space of cycles and $B_{p}\subset Z_{p}$ is the
vector space of boundaries on $M$. The dimension $b_p=\dim H_p(M)$ of the
homology group $H_{p}(M)$ is, by the De Rham theorem, the same Betti number $%
b_p$.

If we know the Betti numbers for all (co)homology groups of the manifold $M$%
, we can calculate the \textit{Euler--Poincar\'{e} characteristic} of $M$ as
\begin{equation*}
\chi (M)=\sum_{p=1}^{n}(-1)^{p}b_{p}.
\end{equation*}

For example, consider a small portion of the De Rham cochain complex of
Figure \ref{Cohomology} spanning a space-time 4--manifold $M$,
\begin{equation*}
\Omega ^{p-1}(M)\overset{d_{p-1}}{\longrightarrow }\Omega ^{p}(M)\overset{%
d_{p}}{\longrightarrow }\Omega ^{p+1}(M)
\end{equation*}%
As we have seen above, cohomology classifies topological spaces by comparing
two subspaces of ~$\Omega ^{p}$:~ (i) the space of $p-$cocycles, $Z^{p}(M)=%
\func{Ker}d_{p}$, and~ (ii) the space of $p-$coboundaries, $B^{p}(M)=\func{Im%
}d_{p-1}$. Thus, for the cochain complex of any space-time 4--manifold we
have,
\begin{equation*}
d^{2}=0\quad \Rightarrow \quad B^{p}(M)\subset Z^{p}(M),
\end{equation*}%
that is, every $p-$coboundary is a $p-$cocycle. Whether the converse of this
statement is true, according to Poincar\'{e} lemma, depends on the
particular topology of a space-time 4--manifold. If every $p-$cocycle is a $%
p-$coboundary, so that $B^{p}$ and $Z^{p}$ are equal, then the cochain
complex is exact at $\Omega ^{p}(M)$. In topologically interesting regions
of a space-time manifold $M$, exactness may fail \cite{Wise}, and we measure
the failure of exactness by taking the $p$th cohomology group
\begin{equation*}
H^{p}(M)=Z^{p}(M)/B^{p}(M).
\end{equation*}

\subsection{Hodge star operator}

The \textit{Hodge star} operator $\star :\Omega ^{p}(M)\rightarrow \Omega
^{n-p}(M)$, which maps any $p-$form $\alpha$ into its \emph{dual} $(n-p)-$%
form $\star\,\alpha$ on a smooth $n-$manifold $M$, is defined as (see, e.g.
\cite{Hodge})
\begin{eqnarray*}
&&\alpha \wedge \star \,\beta =\beta \wedge \star \,\alpha =\left\langle
\alpha ,\beta \right\rangle \mu ,\qquad \text{for }\alpha ,\beta \in \Omega
^{p}(M), \\
&&\star\star\alpha=(-1)^{p(n-p)}\alpha, \\
&&\star\,(c_1\alpha + c_2\beta) = c_1(\star\,\alpha) + c_2(\star\,\beta), \\
&& \alpha\wedge\star\,\alpha = 0 \Rightarrow \alpha\equiv 0.
\end{eqnarray*}
The $\star $ operator depends on the Riemannian metric $g=g_{ij}$ on $M$%
\footnote{%
In local coordinates on a smooth manifold $M$, the metric $g=g_{ij}$ is
defined for any orthonormal basis $(\partial_i=\partial_{x^i})$ in $M$ by
\begin{equation*}
g_{ij}=g(\partial_{i},\partial_{j})=\delta _{ij}, \qquad
\partial_{k}g_{ij}=0.
\end{equation*}%
} and also on the orientation (reversing orientation will change the sign)%
\footnote{%
Hodge $\star$ operator is defined locally in an orthonormal basis (coframe)
of 1--forms $e_i dx^i$ on a smooth manifold $M$ as:
\begin{equation*}
\star(e_i\wedge e_j) = e_k,\qquad (\star)^2 = 1.
\end{equation*}%
} \cite{GaneshADG}. The \textit{volume form} $\mu $ is defined in local
coordinates on an $n-$manifold $M$ as (compare with Hodge inner product
below)
\begin{equation}
\mu =\text{vol}=\star\,(1)=\sqrt{\det (g_{ij})}~dx^{1}\wedge ...\wedge
dx^{n},  \label{mu}
\end{equation}
and the total volume on $M$ is given by
\begin{equation*}
\limfunc{vol}(M)=\int_{M}\star\,(1).
\end{equation*}

For example, in Euclidean $\mathbb{R}^3$ space with Cartesian $(x,y,z)$
coordinates, we have:
\begin{equation*}
\star\, dx=dy\wedge dz,\qquad \star\, dy=dz\wedge dx,\qquad \star\,
dz=dx\wedge dy.
\end{equation*}
The Hodge dual in this case clearly corresponds to the 3D cross--product.

In the 4D--electrodynamics, the dual 2--form \emph{Maxwell} $\star\,F$
satisfies the electric Maxwell equation with the source \cite{MTW},
\begin{equation*}
\text{Dual Bianchi identity}:\quad{d}\,{\star\, F}=\star\,J,
\end{equation*}%
where ${\star \,J}$ is the 3--form dual to the charge--current 1--form ${J}$.

\subsection{Hodge inner product}

For any two $p-$forms $\alpha ,\beta \in \Omega ^{p}(M)$ with compact
support on an $n-$manifold $M$, we define bilinear and positive--definite
Hodge $L^{2}-$inner product as
\begin{equation}
(\alpha ,\beta )=\int_{M}\langle \alpha ,\beta \rangle \star
(1)=\int_{M}\alpha \wedge \star \,\beta ,  \label{L2}
\end{equation}%
where $\alpha \wedge \star \,\beta $ is an $n-$form. We can extend the
product $(\cdot ,\cdot )$ to $L^{2}(\Omega ^{p}(M))$; it remains bilinear
and positive--definite, because as usual in the definition of $L^{2}$,
functions that differ only on a set of measure zero are identified. The
inner product (\ref{L2}) is evidently linear in each variable and symmetric,
$(\alpha ,\beta )=(\beta ,\alpha )$. We have: $(\alpha ,\alpha )\geq 0$ and $%
(\alpha ,\alpha )=0$ iff $\alpha =0$. Also, $(\star\,\alpha,\star\,\beta)=(%
\alpha,\beta)$. Thus, operation (\ref{L2}) turns the space $\Omega ^{p}(M)$
into an infinite--dimensional inner--product space.

From (\ref{L2}) it follows that for every $p-$form $\alpha \in \Omega ^{p}(M)
$ we can define the \emph{norm functional}
\begin{equation*}
\left\Vert \alpha \right\Vert ^{2}=(\alpha ,\alpha )=\int_{M}\langle \alpha
,\alpha \rangle \star (1)=\int_{M}\alpha \wedge \star \,\alpha ,
\end{equation*}%
for which the \emph{Euler--Lagrangian equation} becomes the Laplace equation
(see Hodge Laplacian below),
\begin{equation*}
\Delta \alpha =0.
\end{equation*}

For example, the standard \textit{Lagrangian} for the free Maxwell
electromagnetic field,\newline
$F=dA$ (where $A=A_idx^i$ is the electromagnetic potential 1--form), is
given by \cite{GaneshSprBig,GaneshADG,QuLeap}
\begin{equation*}
\mathcal{L}(A)=\frac{1}{2}(F\wedge \star \,F),
\end{equation*}%
with the corresponding action%
\begin{equation*}
S(A)=\frac{1}{2}\int F\wedge \star \,F.
\end{equation*}%
Using the Hodge $L^{2}-$inner product (\ref{L2}), we can rewrite this
electrodynamic action as%
\begin{equation}
S(A)=\frac{1}{2}(F,F).  \label{act}
\end{equation}

\subsection{Hodge codifferential operator}

The Hodge dual (or, formal adjoint) to the exterior derivative $d:\Omega
^{p}(M)\rightarrow \Omega ^{p+1}(M)$ on a smooth manifold $M$ is the \emph{%
codifferential} $\delta $, a linear map $\delta :\Omega ^{p}(M)\rightarrow
\Omega ^{p-1}(M)$, which is a generalization of the divergence, defined by
\cite{De Rham,Hodge}
\begin{equation*}
\delta =(-1)^{n(p+1)+1}\star d\,\star ,\qquad \text{so that\qquad }%
d=(-1)^{np}\star \delta \star .
\end{equation*}%
That is, if the dimension $n$ of the manifold $M$ is even, then $%
\delta=-\star d \,\star$.

Applied to any $p-$form $\omega \in \Omega ^{p}(M)$, the codifferential $%
\delta $ gives
\begin{equation*}
\delta \omega =(-1)^{n(p+1)+1}\star d\star \omega ,\qquad \delta d\omega
=(-1)^{np+1}\star d\star d\omega .
\end{equation*}%
If $\omega =f$ is a $0-$form, or function, then $\delta f=0$. If a $p-$form $%
\alpha$ is a codifferential of a $(p+1)-$form $\beta$, that is $\alpha
=\delta \beta$, then $\beta$ is called the \emph{coexact} form. A $p-$form $%
\alpha$ is \emph{coclosed} if $\delta\alpha=0$; then $\star\,\alpha$ is
closed (i.e., $d\star\alpha=0$) and conversely.

The Hodge codifferential $\delta $ satisfies the following set of rules:

\begin{itemize}
\item $\delta \delta =\delta ^{2}=0,$ \ the same as $\ dd=d^{2}=0;$

\item $\delta \star =(-1)^{p+1}\star d$; \ $\star \,\delta =(-1)^{p}\star d$;

\item $d\delta \star =\star \,\delta d$; \ $\star \,d\delta =\delta d\star$.
\end{itemize}

Standard example is classical electrodynamics, in which the \textit{gauge
field} is an electromagnetic potential 1--form (a \textit{connection} on a $%
U(1)-$bundle),
\begin{equation*}
A=A_{\mu }dx^{\mu }=A_{\mu }dx^{\mu }+df,\qquad (f=~\text{arbitrary scalar
field}),
\end{equation*}%
with the corresponding electromagnetic field 2--form (the \textit{curvature}
of the connection $A$)
\begin{eqnarray*}
F &=&dA,\qquad \text{in components given by} \\
F &=&\frac{1}{2}F_{\mu \nu }\,dx^{\mu }\wedge dx^{\nu },\qquad \text{with \
\ \ }F_{\mu \nu }=\partial _{\nu }A_{\mu }-\partial _{\mu }A_{\nu }.
\end{eqnarray*}
Electrodynamics is governed by the \emph{Maxwell equations,} which in
exterior formulation read\footnote{%
The first, sourceless Maxwell equation, $dF=0$, gives vector magnetostatics
and magnetodynamics,
\begin{eqnarray*}
\text{Magnetic Gauss' law} &:&\func{div}\mathbf{B}=0,\qquad \\
\text{Faraday's law} &\text{:}&\partial _{t}\mathbf{B}+\func{curl}\mathbf{E}%
=0.
\end{eqnarray*}%
The second Maxwell equation with source, $\delta F=J$ (or, $d\star F = -
\star J$), gives vector electrostatics and electrodynamics,
\begin{eqnarray*}
\text{Electric Gauss' law} &:&\func{div}\mathbf{E}=4\pi \rho ,\qquad \\
\text{Amp\`{e}re's law} &:&\partial _{t}\mathbf{E}-\func{curl}\mathbf{B}%
=-4\pi \mathbf{j}.
\end{eqnarray*}%
}
\begin{eqnarray*}
dF &=&0,\qquad \delta F=-4\pi J,\qquad \text{or in components,} \\
F_{[\mu \nu ,\eta ]} &=&0,\qquad F_{\mu \nu },^{\mu }=-4\pi J_{\mu },
\end{eqnarray*}%
where comma denotes the partial derivative and the 1--form of electric
current $J=J_{\mu }dx^{\mu }$ is conserved, by the electrical \textit{%
continuity equation},
\begin{equation*}
\delta J=0,\qquad \text{or in components,\qquad }J_{\mu },^{\mu }=0.
\end{equation*}

\subsection{Hodge Laplacian operator}

The codifferential $\delta $ can be coupled with the exterior derivative $d$
to construct the \emph{Hodge Laplacian} $\Delta :$ $\Omega
^{p}(M)\rightarrow \Omega ^{p}(M),$ a harmonic generalization of the
Laplace--Beltrami differential operator, given by\footnote{%
Note that the difference $d-\delta ={\partial _{D}}$ is called the \textit{%
Dirac operator}. Its square $\partial _{D}^{2}$\ equals the Hodge Laplacian $%
\Delta $.
\par
Also, in his QFT--based rewriting the Morse topology, E. Witten~\cite%
{WittenMorse} considered also the operators:
\begin{eqnarray*}
d_{t}\;=\;\mathrm{e}^{-tf}de^{tf},\qquad \text{their adjoints}:\quad
d_{t}^{\ast }&=&\mathrm{e}^{tf}de^{-tf}, \\
\text{as well as their Laplacian:\quad }\Delta _{t}&=&d_{t}d_{t}^{\ast
}+d_{t}^{\ast }d_{t}.
\end{eqnarray*}
For $t=0$, $\Delta _{0}$ is the \emph{Hodge Laplacian}, whereas for $%
t\rightarrow \infty $, one has the following expansion
\begin{equation*}
\Delta _{t}=dd^{\ast }+d^{\ast }d+t^{2}\left\Vert df\right\Vert
^{2}+t\sum_{k,j}\frac{\partial ^{2}h}{\partial x^{k}\partial x^{j}}%
[i\,\partial_{x^{k}},dx^{j}],
\end{equation*}%
where $(\partial_{x^{k}})_{k=1,...,n}$ is an orthonormal frame at the point
under consideration. This becomes very large for $t\rightarrow \infty $,
except at the critical points of $f$, i.e., where $df=0$. Therefore, the
eigenvalues of $\Delta _{t}$ will concentrate near the critical points of $f$
for $t\rightarrow \infty $, and we get an \emph{interpolation} between De
Rham cohomology and Morse cohomology.}
\begin{equation*}
\Delta =\delta d+d\delta = (d + \delta)^2.
\end{equation*}%
$\Delta $ satisfies the following set of rules:
\begin{equation*}
\delta \,\Delta =\Delta \,\delta = \delta d\delta;\qquad d\,\Delta =\Delta
\,d = d\delta d;\qquad \star \;\Delta =\Delta \star .
\end{equation*}

A $p-$form $\alpha $ is called \textit{harmonic} iff \
\begin{equation*}
\Delta \alpha =0~\Leftrightarrow ~(d\alpha =0,\delta \alpha =0).
\end{equation*}
Thus, $\alpha $ is harmonic in a compact domain $D\subset M$\footnote{%
A domain $D$ is compact if every open cover of $D$ has a finite subcover.}
iff it is both closed and coclosed in $D$. Informally, every harmonic form
is both closed and coclosed. As a proof, we have:
\begin{equation*}
0=\left( \alpha ,~\Delta \alpha \right) =\left( \alpha ,d\delta \alpha
\right) +\left( \alpha ,\delta d\alpha \right) =\left( \delta \alpha ,\delta
\alpha \right) +\left( d\alpha ,d\alpha \right) .
\end{equation*}%
Since $\left( \beta ,~\beta \right) \geq 0$ for any form $\beta $, $\left(
\delta \alpha ,\delta \alpha \right) $ and $\left( d\alpha ,d\alpha \right) $
must vanish separately. Thus, $d\alpha =0$ and $\delta \alpha =0.$

All harmonic $p-$forms on a smooth manifold $M$ form the vector space $%
H^p_\Delta(M)$.

Also, given a $p-$form $\lambda ,$ there is another $p-$form $\eta $ such
that the equation%
\begin{equation*}
\Delta \eta =\lambda
\end{equation*}%
is satisfied iff for any harmonic $p-$form $\gamma $ we have $\ (\gamma
,\lambda )=0.$

For example, to translate notions from standard 3D vector calculus, we first
identify scalar functions with 0--forms, field intensity vectors with
1--forms, flux vectors with 2--forms and scalar densities with 3--forms. We
then have the following correspondence:

grad $\longrightarrow $ $d$ : \ on 0--forms;\qquad \qquad \qquad curl $%
\longrightarrow $ $\star \,d$ : \ on 1--forms;

div $\longrightarrow $ $\delta $ : \ on 1--forms;\qquad \qquad \qquad\ \ div
grad $\longrightarrow $ $\Delta $ : \ on 0--forms;

\qquad \qquad \qquad curl curl $-$ grad div $\longrightarrow \Delta $ : \ on
1--forms.\newline

We remark here that exact and coexact $p-$forms ($\alpha=d\beta$ and $%
\omega=\delta\beta$) are mutually orthogonal with respect to the $L^{2}-$%
inner product (\ref{L2}). The orthogonal complement consists of forms that
are both closed and coclosed: that is, of harmonic forms ($\Delta\gamma=0$).

\subsection{Hodge adjoints and self--adjoints}

If $\alpha $ is a $p-$form and $\beta $ is a $(p+1)-$form then we have \cite%
{De Rham}
\begin{equation}
(d\alpha ,\beta )=(\alpha ,\delta \beta )\qquad\mathrm{and}\qquad
(\delta\alpha,\beta)=(\alpha,d\beta).  \label{adj}
\end{equation}%
This relation is usually interpreted as saying that the two exterior
differentials, $d$ and $\delta \,,$ are \textit{adjoint} (or, dual) to each
other. This identity follows from the fact that for the volume form $\mu $
given by (\ref{mu}) we have $d\mu=0$ and thus
\begin{equation*}
\int_M d(\alpha \wedge \star\beta)=0.
\end{equation*}

Relation (\ref{adj}) also implies that the Hodge Laplacian $\Delta $ is
\textit{self--adjoint} (or, self--dual),
\begin{equation*}
(\Delta \alpha ,\beta )=(\alpha ,\Delta \beta ),
\end{equation*}%
which is obvious as either side is $(d\alpha ,d\beta )+(\delta \alpha
,\delta \beta ).$ Since $(\Delta \alpha ,\alpha )\geq 0,$ with $(\Delta
\alpha ,\alpha )=0$ only when $\Delta \alpha =0,$ $\Delta $ is a
positive--definite (elliptic) self--adjoint differential operator.

\subsection{Hodge decomposition theorem}

The celebrated \emph{Hodge decomposition theorem} (HDT) states that, on a
compact orientable smooth $n-$manifold $M$ (with $n\geq p$), any exterior $%
p- $form can be written as a unique sum of an \emph{exact} form, a \emph{%
coexact} form, and a \emph{harmonic} form. More precisely, for any form $%
\omega \in \Omega ^{p}(M)$ there are unique forms $\alpha \in \Omega
^{p-1}(M),$ $\beta \in \Omega ^{p+1}(M)$ and a harmonic form $\gamma \in
\Omega ^{p}(M),$ such that \bigbreak%
\centerline{\fbox{\parbox{11cm}{
{\large\begin{equation*} {\rm HDT:}\qquad\stackrel{\rm
any\,form}{\omega }\quad =\quad
\stackrel{\rm exact}{d\alpha }\quad + \quad\stackrel{\rm coexact}{\delta\beta }\quad + \quad\stackrel{\rm harmonic}{\gamma }
\end{equation*}}}}}\bigbreak\noindent For the proof, see \cite{Hodge,De Rham}%
.

In physics community, the exact form $d\alpha$ is called \emph{longitudinal}%
, while the coexact form $\delta\beta$ is called \emph{transversal}, so that
they are mutually orthogonal. Thus any form can be orthogonally decomposed
into a harmonic, a longitudinal and transversal form. For example, in fluid
dynamics, any vector-field $v$ can be decomposed into the sum of two
vector-fields, one of which is divergence--free, and the other is curl--free.

Since $\gamma $ is harmonic, $d\gamma =0.$ Also, by Poincar\'{e} lemma, $%
d(d\alpha )=0.$ In case $\omega $ is a closed $p-$form, $d\omega =0,$ then
the term $\delta \beta $ in HDT is absent, so we have the \emph{short Hodge
decomposition},
\begin{equation}
\omega =d\alpha +\gamma ,  \label{sHd}
\end{equation}%
thus $\omega $ and $\gamma $ differ by $d\alpha$. In topological
terminology, $\omega $ and $\gamma$ belong to the same \textit{cohomology
class} $[\omega ]\in H^{p}(M)$. Now, by the De Rham theorems it follows that
if $C$ is any $p-$cycle, then
\begin{equation*}
\int_{C}\omega =\int_{C}\gamma ,
\end{equation*}%
that is, $\gamma $ and $\omega $ have the same periods$.$ More precisely, if
$\omega $ is any closed $p-$form, then there exists a unique harmonic $p-$%
form $\gamma $ with the same periods as those of $\omega $ (see \cite{De
Rham,Flanders}).

The \emph{Hodge--Weyl theorem} \cite{Hodge,De Rham} states that every De
Rham cohomology class has a unique harmonic representative. In other words,
the space $H^p_\Delta(M)$ of harmonic $p-$forms on a smooth manifold $M$ is
isomorphic to the De Rham cohomology group (\ref{DR}), or\newline
$H^p_\Delta(M) \cong H^p_{DR}(M)$. That is, the harmonic part $\gamma$ of
HDT depends only on the global structure, i.e., the topology of $M$.

For example, in $(2+1)$D electrodynamics, $p-$form Maxwell equations in the
Fourier domain $\Sigma $ are written as \cite{Teixeira}%
\begin{equation*}
\begin{array}{ll}
dE=i\omega B, & \qquad dB=0, \\
dH=-i\omega D+J, & \qquad dD=Q,%
\end{array}%
\end{equation*}%
where $H$ is a 0--form (magnetizing field), $D$ (electric displacement
field), $J$ (electric current density) and $E$ (electric field) are
1--forms, while $B$ (magnetic field) and $Q$ (electric charge density) are
2--forms. From $d^{2}=0$ it follows that the $J$ and the $Q$ satisfy the
\textit{continuity equation}%
\begin{equation*}
dJ=i w Q,
\end{equation*}
where $i=\sqrt{-1}$ and $w$ is the field frequency. Constitutive equations,
which include all metric information in this framework, are written in terms
of Hodge star operators (that fix an isomorphism between $p$ forms and $%
(2-p) $ forms in the $(2+1)$ case)%
\begin{equation*}
D=\star \,E,\qquad B=\star \,H.
\end{equation*}

Applying HDT to the electric field intensity 1--form $E$, we get \cite{He}%
\begin{equation*}
E=d\phi +\delta A+\chi ,
\end{equation*}%
where $\phi $ is a 0--form (a scalar field) and $A$ is a 2--form; $d\phi $
represents the static field and $\delta A$ represents the dynamic field, and
$\chi $ represents the harmonic field component. If domain $\Sigma $ is
contractible, $\chi $ is identically zero and we have the short Hodge
decomposition,%
\begin{equation*}
E=d\phi +\delta A.
\end{equation*}

\section{Hodge Decomposition and Gauge Path Integral}

\subsection{Feynman Path Integral}

The `driving engine' of quantum field theory is the Feynman path integral.
Very briefly, there are three basic forms of the path integral (see, e.g.,
\cite{QuLeap}):

1. \emph{Sum--over--histories,} developed in Feynman's version of quantum
mechanics (QM)\footnote{%
Feynman's \emph{amplitude} is a space-time version of the Schr\"{o}dinger's
\emph{wave-function} $\psi $, which describes how the (non-relativistic)
quantum state of a physical system changes in space and time, i.e.,
\begin{equation*}
\langle \mathrm{Out}_{t_{fin}}|\mathrm{In}_{t_{ini}}\rangle =\psi (\mathbf{x}%
,t),\qquad (\text{for~~}\mathbf{x}\in \lbrack \mathrm{In},\mathrm{Out}%
],~t\in \lbrack t_{ini},t_{fin}]).
\end{equation*}%
In particular, quantum wave-function $\psi $ is a complex--valued function
of real space variables $\mathbf{x}=(x_{1},x_{2},...,x_{n})\in \mathbb{R}^{n}
$, which means that its domain is in $\mathbb{R}^{n}$ and its range is in
the complex plane, formally $\psi (\mathbf{x}):\mathbb{R}^{n}\rightarrow
\mathbb{C}.$ For example, the one--dimensional \emph{stationary plane wave}
with wave number $k$ is defined as
\begin{equation*}
\psi (x)=\mathrm{e}^{\mathrm{i}kx},\qquad (\mathrm{for~~}x\in \mathbb{R}),
\end{equation*}%
where the real number $k$ describes the wavelength, $\lambda =2\pi /k.$ In $%
n $ dimensions, this becomes
\begin{equation*}
\psi (x)=\mathrm{e}^{\mathrm{i}\mathbf{p\cdot x}},
\end{equation*}%
where the momentum vector $\mathbf{p=k}$ is the vector of the wave numbers $%
\mathbf{k}$ in natural units (in which $\hbar =m=1$).
\par
More generally, quantum wave-function is also time dependent, $\psi =\psi (%
\mathbf{x,t}).$ The time--dependent plane wave is defined by%
\begin{equation}
\psi (\mathbf{x,t})=\mathrm{e}^{\mathrm{i}\mathbf{p\cdot x-}\mathrm{i}%
p^{2}t/2}.  \label{plW}
\end{equation}%
\par
In general, $\psi (\mathbf{x,t})$ is governed by the Schr\"{o}dinger
equation \cite{VQM,QuLeap} (in natural units $\hbar =m=0$)
\begin{equation}
\mathrm{i}\frac{\partial }{\partial t}\psi (\mathbf{x},t)=-\frac{1}{2}\Delta
\psi (\mathbf{x},t),  \label{schgen}
\end{equation}%
where $\Delta $ is the $n-$dimensional Laplacian. The solution of (\ref%
{schgen}) is given by the integral of the time--dependent plane wave (\ref%
{plW}),%
\begin{equation*}
\psi (\mathbf{x},t)=\frac{1}{(2\pi )^{n/2}}\int_{\mathbb{R}^{n}}\mathrm{e}^{%
\mathrm{i}\mathbf{p\cdot x-}\mathrm{i}p^{2}t/2}\hat{\psi}_{0}(\mathbf{p}%
)d^{n}p,
\end{equation*}%
which means that $\psi (\mathbf{x},t)$ is the inverse Fourier transform of
the function%
\begin{equation*}
\hat{\psi}(\mathbf{p},t)=\mathrm{e}^{\mathbf{-}\mathrm{i}p^{2}t/2}\hat{\psi}%
_{0}(\mathbf{p}),
\end{equation*}%
where $\hat{\psi}_{0}(\mathbf{p})$ has to be calculated for each initial
wave-function. For example, if initial wave-function is Gaussian,
\begin{equation*}
f(x)=\exp (-a\frac{x^{2}}{2}),\qquad \text{with the Fourier transform}\qquad
\hat{f}(p)=\frac{1}{\sqrt{a}}\exp (-\frac{p^{2}}{2a}).
\end{equation*}%
$\hspace{2cm}\text{ then}\qquad \hat{\psi}_{0}(p)=\frac{1}{\sqrt{a}}\exp (-%
\frac{p^{2}}{2a}).$} \cite{FeynQM};

2. \emph{Sum--over--fields,} started in Feynman's version of quantum
electrodynamics (QED) \cite{FeynQED} and later improved by Fadeev--Popov \cite{Faddeev};

3. \emph{Sum--over--geometries/topologies} in quantum gravity (QG),
initiated by S. Hawking and properly developed in the form of causal
dynamical triangulations (see \cite{Ambjorn}; for a `softer' review, see
\cite{Loll}).

In all three versions, Feynman's \emph{action--amplitude formalism} includes
two components:

1. A real--valued, classical, \emph{Hamilton's action functional,}
\begin{equation*}
S[\Phi ]~:=~\int_{t_{ini}}^{t_{fin}}{L}[\Phi ]\,dt,
\end{equation*}%
with the Lagrangian energy function defined over the Lagrangian density $%
\mathcal{L}$,
\begin{equation*}
{L}[\Phi ]=\int d^{n}x\,\mathcal{L}(\Phi ,\partial _{\mu }\Phi ),\qquad
(\partial _{\mu }\equiv \partial /\partial x^{\mu }),
\end{equation*}%
while $\Phi $ is a common symbol denoting all three things to be summed upon
(histories, fields and geometries). The action functional $S[\Phi ]$ obeys
the \emph{Hamilton's least action principle,} ~$\delta S[\Phi ]=0,$~ and
gives, using standard variational methods,\footnote{%
In Lagrangian field theory, the fundamental quantity is the action
\begin{equation*}
S[\Phi ]=\int_{t_{in}}^{t_{out}}L\,dt=\int_{\mathbb{R}^{4}}d^{n}x\,\mathcal{L%
}(\Phi ,\partial _{\mu }\Phi )\,,
\end{equation*}%
so that the least action principle, $\delta S[\Phi ]=0,$ gives
\begin{eqnarray*}
0 &=&\int_{\mathbb{R}^{4}}d^{n}x\,\left\{ \frac{\partial \mathcal{L}}{%
\partial \Phi }\delta \Phi +\frac{\partial \mathcal{L}}{\partial (\partial
_{\mu }\Phi )}\delta (\partial _{\mu }\Phi )\right\} \\
&=&\,\int_{\mathbb{R}^{4}}d^{n}x\,\left\{ \frac{\partial \mathcal{L}}{%
\partial \Phi }\delta \Phi -\partial _{\mu }\left( \frac{\partial \mathcal{L}%
}{\partial (\partial _{\mu }\Phi )}\right) \delta \Phi +\partial _{\mu
}\left( \frac{\partial \mathcal{L}}{\partial (\partial _{\mu }\Phi )}\delta
\Phi \right) \right\} .
\end{eqnarray*}%
The last term can be turned into a surface integral over the boundary of the
$\mathbb{R}^{4}$ (4D space-time region of integration). Since the initial
and final field configurations are assumed given, $\delta \Phi =0$ at the
temporal beginning $t_{in}$ and end $t_{out}$\ of this region, which implies
that the surface term is zero. Factoring out the $\delta \Phi $ from the
first two terms, and since the integral must vanish for arbitrary $\delta
\Phi $, we arrive at the Euler-lagrange equation of motion for a field,
\begin{equation*}
\partial _{\mu }\left( \frac{\partial \mathcal{L}}{\partial (\partial _{\mu
}\Phi )}\right) -\frac{\partial \mathcal{L}}{\partial \Phi }=0.
\end{equation*}%
If the Lagrangian (density) $\mathcal{L}$\ contains more fields, there is
one such equation for each. The momentum density $\pi (x)$ of a field,
conjugate to $\Phi (x)$ is defined as: ~ $\pi (x)=\frac{\partial \mathcal{L}%
}{\partial _{\mu }\Phi (x)}.$%
\par
For example, the standard electromagnetic action
\begin{equation*}
S=-\frac{1}{4}\int_{\mathbb{R}^{4}}d^{4}x\,F_{\mu \nu }F^{\mu \nu },\qquad
\text{where}\qquad F_{\mu \nu }=\partial _{\mu }A_{\nu }-\partial _{\nu
}A_{\mu },
\end{equation*}%
gives the sourceless Maxwell's equations: ~
\begin{equation*}
\partial _{\mu }F^{\mu \nu }=0,\qquad \epsilon ^{\mu \nu \sigma \eta
}\partial _{\nu }F_{\sigma \eta }=0,
\end{equation*}%
where the field strength tensor $F_{\mu \nu }$ and the Maxwell equations are
invariant under the \emph{gauge transformations,}
\begin{equation*}
A_{\mu }\longrightarrow A_{\mu }+\partial _{\mu }\epsilon .
\end{equation*}%
\par
The equations of motion of charged particles are given by the Lorentz--force
equation,
\begin{equation*}
m{\frac{du^{\mu }}{d\tau }}=eF^{\mu \nu }u_{\nu },
\end{equation*}%
where $e$ is the charge of the particle and $u^{\mu }(\tau )$ its
four-velocity as a function of the proper time.} the Euler--Lagrangian
equations, which define the shortest path, the extreme field, and the
geometry of minimal curvature (and without holes).

2. A complex--valued, quantum \emph{transition amplitude},\footnote{%
The transition amplitude is closely related to \emph{partition function} $Z,$
which is a quantity that encodes the statistical properties of a system in
thermodynamic equilibrium. It is a function of temperature and other
parameters, such as the volume enclosing a gas. Other thermodynamic
variables of the system, such as the total energy, free energy, entropy, and
pressure, can be expressed in terms of the partition function or its
derivatives. In particular, the partition function of a \emph{canonical
ensemble}~is defined as a sum ~~$Z(\beta )=\sum_{j}\mathrm{e}^{-\beta
E_{j}}, $~ where $\beta =1/(k_{B}T)$ is the `inverse temperature', where $T$
is an ordinary temperature and $k_{B}$ is the Boltzmann's constant. However,
as the position $x^{i}$ and momentum $p_{i}$ variables of an $i$th particle
in a system can vary continuously, the set of microstates is actually
uncountable. In this case, some form of \textit{coarse--graining} procedure
must be carried out, which essentially amounts to treating two mechanical
states as the same microstate if the differences in their position and
momentum variables are `small enough'. The partition function then takes the
form of an integral. For instance, the partition function of a gas
consisting of $N$ molecules is proportional to the $6N-$dimensional
phase--space integral,
\begin{equation*}
Z(\beta )\sim \int_{\mathbb{R}^{6N}}\,d^{3}p_{i}\,d^{3}x^{i}\exp [-\beta
H(p_{i},x^{i})],
\end{equation*}%
where $H=H(p_{i},x^{i}),$ ($i=1,...,N$) \ is the classical Hamiltonian
(total energy) function.
\par
Given a set of random variables $X_{i}$ taking on values $x^{i}$, and purely
potential Hamiltonian function $H(x^{i})$, the partition function is defined
as
\begin{equation*}
Z(\beta )=\sum_{x^{i}}\exp \left[ -\beta H(x^{i})\right] .
\end{equation*}%
The function $H$ is understood to be a real-valued function on the space of
states $\{X_{1},X_{2},\cdots \}$ while $\beta $ is a real-valued free
parameter (conventionally, the inverse temperature). The sum over the $x^{i}$
is understood to be a sum over all possible values that the random variable $%
X_{i}$ may take. Thus, the sum is to be replaced by an integral when the $%
X_{i}$ are continuous, rather than discrete. Thus, one writes
\begin{equation*}
Z(\beta )=\int dx^{i}\exp \left[ -\beta H(x^{i})\right] ,
\end{equation*}%
for the case of continuously-varying random variables $X_{i}$.
\par
Now, the number of variables $X_{i}$ need not be countable, in which case
the set of coordinates $\{x^{i}\}$ becomes a field\ $\phi =\phi (x),$ so\
the sum is to be replaced by the \emph{Euclidean path integral} (that is a
Wick--rotated Feynman transition amplitude (\ref{Wick}) in imaginary time),
as
\begin{equation*}
Z(\phi )=\int \mathcal{D}[\phi ]\exp \left[ -H(\phi )\right] .
\end{equation*}%
\par
More generally, in quantum field theory, instead of the field Hamiltonian $%
H(\phi )$ we have the action $S(\phi )$ of the theory. Both Euclidean path
integral,
\begin{equation}
Z(\phi )=\int \mathcal{D}[\phi ]\exp \left[ -S(\phi )\right] ,\qquad \text{%
real path integral in imaginary time,}  \label{Eucl}
\end{equation}%
and Lorentzian one,
\begin{equation}
Z(\phi )=\int \mathcal{D}[\phi ]\exp \left[ \mathrm{i}S(\phi )\right]
,\qquad \text{complex path integral in real time,}  \label{Lor}
\end{equation}%
are usually called `partition functions'. While the Lorentzian path integral
(\ref{Lor}) represents a quantum-field theory-generalization of the Schr\"{o}%
dinger equation, the Euclidean path integral (\ref{Eucl}) represents a
statistical-field-theory generalization of the Fokker--Planck equation.}
\begin{equation}
\langle \mathrm{Out}_{t_{fin}}|\mathrm{In}_{t_{ini}}\rangle ~:=\int_{\mathrm{%
\Omega }}\mathcal{D}[\Phi ]\,\mathrm{e}^{\mathrm{i}S[\Phi ]},
\label{pathInt1}
\end{equation}%
where $\mathcal{D}[\Phi ]$ is `an appropriate' Lebesgue--type measure,
\begin{equation*}
\mathcal{D}[\Phi ]=\lim_{N\rightarrow \infty }\prod_{s=1}^{N}\Phi
_{s}^{i},\qquad (i=1,...,n),
\end{equation*}%
so that we can `safely integrate over a continuous spectrum and sum over a
discrete spectrum of our problem domain $\Omega $', of which the absolute
square is the real--valued probability density function,
\begin{equation*}
P~:=~|\langle \mathrm{Out}_{t_{fin}}|\mathrm{In}_{t_{ini}}\rangle \rangle
|^{2}.
\end{equation*}

This procedure can be redefined in a mathematically cleaner way if we
Wick--rotate the time variable $t$ to imaginary values, $t\mapsto \tau={} t$%
, thereby making all integrals real:
\begin{equation}
\int \mathcal{D}[\Phi]\, \mathrm{e}^{\mathrm{i} S[\Phi]}~\cone{Wick}\quad
\int \mathcal{D}[\Phi]\, \mathrm{e}^{-S[\Phi]}.  \label{Wick}
\end{equation}

For example, in non-relativistic quantum mechanics (see Appendix), the
propagation amplitude from $x_{a}$\ to $x_{b}$\ is given by the \emph{%
configuration path integral}\footnote{%
On the other hand, the \emph{phase--space path integral} (without peculiar
constants in the functional measure) reads
\begin{equation*}
U(q_{a},q_{b};T)=\left( \prod_{i}\int \mathcal{D}[q(t)]\mathcal{D}%
[p(t)]\right) \exp \left[ \mathrm{i}\int_{0}^{T}\left( p_{i}\dot{q}%
^{i}-H(q,p)\right) \,dt\right] ,
\end{equation*}%
where the functions $q(t)$\ (space coordinates) are constrained at the
endpoints, but the functions $p(t)$ (canonically--conjugated momenta) are
not. The functional measure is just the product of the standard integral
over phase space at each point in time
\begin{equation*}
\mathcal{D}[q(t)]\mathcal{D}[p(t)]=\prod_{i}\frac{1}{2\pi }\int dq^{i}dp_{i}.
\end{equation*}%
Applied to a non-relativistic real scalar field $\phi (x,t)$, this path
integral becomes
\begin{equation*}
\left\langle \phi _{b}(x,t)|\,\mathrm{e}^{-\mathrm{i}HT}|\phi
_{a}(x,t)\right\rangle =\int \mathcal{D}[\phi ]\exp \left[ \mathrm{i}%
\int_{0}^{T}\mathcal{L}(\phi )\,d^{4}x\right] ,\quad \text{with \ \ }%
\mathcal{L}(\phi )=\frac{1}{2}(\partial _{\mu }\phi )^{2}-V(\phi ).
\end{equation*}%
}
\begin{equation*}
U(x_{a},x_{b};T)=\left\langle x_{b}|x_{a}\right\rangle =\left\langle x_{b}|\,%
\mathrm{e}^{-\mathrm{i}HT}|x_{a}\right\rangle =\int \mathcal{D}[x(t)]\,%
\mathrm{e}^{\mathrm{i}S[x(t)]},
\end{equation*}%
which satisfies the Schr\"{o}dinger equation (in natural units)
\begin{equation*}
i\frac{\partial }{\partial T}U(x_{a},x_{b};T)=\hat{H}U(x_{a},x_{b};T),\qquad
\text{where \ \ \ \ \ }\hat{H}=-\frac{1}{2}\frac{\partial ^{2}}{\partial
x_{b}^{2}}+V(x_{b}).
\end{equation*}

\subsubsection{Functional measure on the space of differential forms}

The Hodge inner product (\ref{L2}) leads to a natural (metric--dependent)
functional measure $\mathcal{D}\mu \lbrack \omega ]$ on $\Omega ^{p}(M)$,
which normalizes the \textit{Gaussian functional integral}
\begin{equation}
\int \mathcal{D}\mu \lbrack \omega ]\,\mathrm{e}^{\mathrm{i}\langle \omega
|\omega \rangle }=1.  \label{2.4}
\end{equation}

One can use the invariance of (\ref{2.4}) to determine how the functional
measure transforms under the Hodge decomposition. Using HDT and its
orthogonality with respect to the inner product (\ref{L2}), it was shown in
\cite{Gegenberg} that
\begin{equation}
\langle \omega ,\omega \rangle =\langle \gamma ,\gamma \rangle +\langle
d\alpha ,d\alpha \rangle +\langle \delta \beta ,\delta \beta \rangle
=\langle \gamma ,\gamma \rangle +\langle \alpha ,\delta d\alpha \rangle
+\langle \beta ,d\delta \beta \rangle\, ,  \label{2.5}
\end{equation}
where the following differential/conferential identities were used \cite%
{Choquet}
\begin{equation*}
\langle d\alpha ,d\alpha \rangle =\langle \alpha ,\delta d\alpha \rangle
\qquad \text{and\qquad }\langle \delta \beta ,\delta \beta \rangle =\langle
\beta ,d\delta \beta \rangle .
\end{equation*}
Since, for any linear operator $O$, one has
\begin{equation*}
\int \mathcal{D}\mu \lbrack \omega ]\exp \mathrm{i}\langle \omega |O\omega
\rangle ={\det }^{-1/2}(O),
\end{equation*}
(\ref{2.4}) and (\ref{2.5}) imply that
\begin{equation*}
\mathcal{D}\mu \lbrack \omega ]=\mathcal{D}\mu \lbrack \gamma ]\mathcal{D}%
\mu \lbrack \alpha ]\mathcal{D}\mu \lbrack \beta ]\,{\det }^{1/2}(\delta d){%
\det }^{1/2}(d\delta ).
\end{equation*}

\subsubsection{Abelian Chern--Simons theory}

Recall that the classical action for an Abelian Chern--Simons theory,
\begin{equation*}
S=\int_{M}A\wedge dA\,,
\end{equation*}%
is invariant (up to a total divergence) under the gauge transformation:
\begin{equation}
A\longmapsto A+d\varphi .  \label{2.10}
\end{equation}%
We wish to compute the \textit{partition function} for the theory
\begin{equation*}
Z:=\int \frac{1}{V_{G}}\mathcal{D}\mu \lbrack A]\,\mathrm{e}^{\mathrm{i}%
S[A]},
\end{equation*}%
where $V_{G}$ denotes the volume of the group of gauge transformations in (%
\ref{2.10}), which must be factored out of the partition function in order
to guarantee that the integration is performed only over physically distinct
gauge fields. We can handle this by using the Hodge decomposition to
parametrize the potential $A$ in terms of its gauge invariant, and gauge
dependent parts, so that the volume of the group of gauge transformations
can be explicitly factored out, leaving a functional integral over gauge
invariant modes only \cite{Gegenberg}.

We now transform the integration variables:
\begin{equation*}
A\longmapsto \alpha ,\beta ,\gamma ,
\end{equation*}
where $\alpha ,\beta ,\gamma $ parameterize respectively the exact, coexact,
and harmonic parts of the connection A. Using the Jacobian (\ref{2.5}) as
well as the following identity on 0--forms $\Delta =\delta d,$ we get \cite%
{Gegenberg}
\begin{equation*}
Z=\int \frac{1}{V_{G}}\mathcal{D}\mu \lbrack \alpha ]\mathcal{D}\mu \lbrack
\beta ]\mathcal{D}\mu \lbrack \gamma ]\,{\det }^{1/2}\left( \Delta \right) {%
\det }^{1/2}\left( d\delta \right) \mathrm{e}^{\mathrm{i}S},
\end{equation*}
from which it follows that
\begin{equation}
V_{G}=\int \mathcal{D}\mu \lbrack \alpha ],  \label{2.13}
\end{equation}
while the classical action functional becomes, after integrating by parts,
using the harmonic properties of $\gamma $ and the nilpotency of the
exterior derivative operators, and dropping surface terms:
\begin{equation*}
S=-\langle \beta ,\star \delta d\delta \beta \rangle \,\,.
\end{equation*}
Note that $S$ depends only the coexact (transverse) part of $A$. Using (\ref%
{2.13}) and integrating over $\beta $ yields:
\begin{equation*}
Z=\int \mathcal{D}\mu \lbrack \gamma ]{\det }^{-1/2}\left( \star \delta
d\delta \right) {\det }^{1/2}\left( \Delta \right) {\det }^{1/2}\left(
d\delta \right) .
\end{equation*}
Also, it was proven in \cite{Gegenberg} that
\begin{equation*}
{\det }(\star \delta d\delta )={\det }^{1/2}((d\delta d)(\delta d\delta ))={%
\det }^{\frac{3}{2}}(d\delta ).
\end{equation*}
As a consequence of Hodge duality we have the identity
\begin{equation*}
{\det }(\delta d)={\det }(d\delta ),
\end{equation*}
from which it follows that
\begin{equation*}
Z=\int \mathcal{D}\mu \lbrack \gamma ]\,{\det }^{-3/4}\left( \Delta
_{(1)}^{T}\right) {\det }^{1/2}\left( \Delta \right) {\det }^{1/2}\left(
\Delta _{(1)}^{T}\right) \,.
\end{equation*}
The operator $\Delta _{(1)}^{T}$ is the transverse part of the Hodge
Laplacian acting on $1-$forms:
\begin{equation*}
\Delta _{(1)}^{T}:=(\delta d)_{(1)}.
\end{equation*}
Applying identity for the Hodge Laplacian $\Delta _{(p)}$ \cite{Gegenberg}

\begin{equation*}
\det \left( \Delta _{(p)}\right) =\det \left( (\delta d)_{(p)}\right) \det
\left( (\delta d)_{(p-1)}\right) ,
\end{equation*}
we get
\begin{equation*}
{\det }\left( \Delta _{(1)}^{T}\right) ={\det }\left( \Delta _{(1)}\right) /{%
\det }\left( \Delta \right)
\end{equation*}
and hence
\begin{equation*}
Z=\int \mathcal{D}\mu \lbrack \gamma ]\,{\det }^{-1/4}\left( \Delta
_{(1)}\right) {\det }^{3/4}\left( \Delta \right).  \label{Zfin}
\end{equation*}
The space of harmonic forms $\gamma $\ (of any order) is a finite set.
Hence, the integration over harmonic forms (\ref{Zfin}) is a simple sum.

\subsection{Appendix: Path Integral in Quantum Mechanics}

The \textit{amplitude} $\langle q_{F}|\mathrm{e}^{-\mathrm{i}%
HT}|q_{I}\rangle $\ of a quantum-mechanical system (in natural units $\hbar
=c=1$) with Hamiltonian $H$ to propagate from a point $q_{I}$ to a point $%
q_{F}$ in time $T$ is governed by the \emph{unitary operator} $\mathrm{e}^{-%
\mathrm{i}HT}$, or more completely, by the complex \emph{Lorentzian path
integral} (in real time):
\begin{equation}
\langle q_{F}|\mathrm{e}^{-\mathrm{i}HT}|q_{I}\rangle =\int \mathcal{D}\left[
q(t)\right] \mathrm{e}^{\mathrm{i}\int_{0}^{T}dtL(\dot{q},q)}.  \label{z6}
\end{equation}%
For example, a quantum particle in a potential $V(q)$ has the \emph{%
Hamiltonian operator}:\newline
$\hat{H}=\hat{p}^{2}/2m+V(\hat{q})$, with the corresponding \emph{Lagrangian}%
: $L=\frac{1}{2}m\dot{q}^{2}-V(q)$, so the path integral reads: \
\begin{equation*}
\langle q_{F}|\mathrm{e}^{-\mathrm{i}HT}|q_{I}\rangle =\int \mathcal{D}\left[
q(t)\right] \mathrm{e}^{\mathrm{i}\int_{0}^{T}dt[\frac{1}{2}m\dot{q}%
^{2}-V(q)]}.
\end{equation*}

It is somewhat more rigorous to perform a so-called \textit{Wick rotation}
to Euclidean time, which means substituting $t\rightarrow -\mathrm{i}t$ and
rotating the integration contour in the complex $t-$plane, so that we obtain
the real path integral in complex time:
\begin{equation*}
\langle q_{F}|\mathrm{e}^{-\mathrm{i}HT}|q_{I}\rangle =\int \mathcal{D}\left[
q(t)\right] \mathrm{e}^{\mathrm{-}\int_{0}^{T}dtL(\dot{q},q)}=\int \mathcal{D%
}\left[ q(t)\right] \mathrm{e}^{\mathrm{-}(\mathrm{i}/\hbar )\int_{0}^{T}dt[%
\frac{1}{2}m\dot{q}^{2}-V(q)]},
\end{equation*}%
known as the \textit{Euclidean path integral}.

One particularly nice feature of the path-integral formalism is that the
classical limit of quantum mechanics can be recovered easily. We simply
restore Planck's constant $\hbar $ in (\ref{z6})
\begin{equation*}
\langle q_{F}|\mathrm{e}^{-(\mathrm{i}/\hbar )HT}|q_{I}\rangle =\int
\mathcal{D}\left[ q(t)\right] \mathrm{e}^{(\mathrm{i}/\hbar )\int_{0}^{T}dtL(%
\dot{q},q)},
\end{equation*}%
and take the $\hbar \rightarrow 0$ limit. Applying the stationary phase
method (or steepest descent) , see \cite{ZeeQFT}), we obtain $\mathrm{e}^{(%
\mathrm{i}/\hbar )\int_{0}^{T}dtL(\dot{q}_{c},q_{c})}$,\footnote{%
To do an exponential integral of the form: ~$I=\int_{-\infty }^{+\infty }dq\,%
\mathrm{e}^{\mathrm{-}(\mathrm{i}/\hbar )f(q)},$ we often have to resort to
the following steepest-descent approximation \cite{ZeeQFT}. In the limit of $%
\hbar $ small, this integral is dominated by the minimum of $f(q)$.
Expanding: $f(q)=f(a)+\frac{1}{2}f"(a)(q-a)^{2}+O[(q-a)^{3}],$ and applying
the \emph{Gaussian integral} rule:
\begin{equation*}
\int_{-\infty }^{+\infty }dq\,\mathrm{e}^{\mathrm{-}\frac{1}{2}ax^{2}}=\sqrt{%
\frac{2\pi }{a}},\qquad \text{we obtain:}\qquad I=\mathrm{e}^{\mathrm{-}%
(1/\hbar )f(a)}\left( \frac{2\pi \hbar }{f"(a)}\right) \mathrm{e}^{\mathrm{-}%
O(\hbar ^{\frac{1}{2}})}.
\end{equation*}%
} where $q_{c}(t)$ is the recovered \emph{classical path} determined by
solving the \textit{Euler-Lagrangian equation}: $(d/dt)(\delta L/\delta
q)-(\delta L/\delta q)=0,$ with appropriate boundary conditions.

\end{document}